\documentclass[a4paper,12pt,reqno]{amsart}
\usepackage{amssymb,amsmath}
\usepackage{graphicx}
\usepackage{enumitem}
\usepackage{subfigure}
\usepackage{hyperref}
\usepackage{pgf,tikz}
\usetikzlibrary{arrows}
\usepackage[all]{xy} 
\input{covers.def}
\input{comments.tex} 

\def\g{\sf{g}}
\def\restriction#1#2{\mathchoice
              {\setbox1\hbox{${\displaystyle #1}_{\scriptstyle #2}$}
              \restrictionaux{#1}{#2}}
              {\setbox1\hbox{${\textstyle #1}_{\scriptstyle #2}$}
              \restrictionaux{#1}{#2}}
              {\setbox1\hbox{${\scriptstyle #1}_{\scriptscriptstyle #2}$}
              \restrictionaux{#1}{#2}}
              {\setbox1\hbox{${\scriptscriptstyle #1}_{\scriptscriptstyle #2}$}
              \restrictionaux{#1}{#2}}}
\def\restrictionaux#1#2{{#1\,\smash{\vrule height .8\ht1 depth .85\dp1}}_{\,#2}}
\newcommand{\fonction}[5]{\begin{array}{cccc}
#1: & #2 & \longrightarrow & #3 \\
    & #4 & \longmapsto & #5 \end{array}}
    
\newcommand{\tq}{\, : \,}
\newcommand{\Ms}{M_{\sigma}}
\newcommand{\rie}{Riemann} 
\title{Extension of analytic covers by extension of their ramification divisors}
\author{L. Lavoine}
\date{\today}

\address{
Universit\'e de Lille-1, UFR de Math\'ematiques, 59655 Villeneuve
d'Ascq, France} \email{landry.lavoine@math.univ-lille1.fr}


\begin{document}
\begin{abstract}
This paper deals with extension of analytic covers.
We prove topological extension theorems for analytic covers. The main result is an extension theorem which only uses the extension of the ramification divisor. We give also a Thullen-type and a Hartogs-type extensions theorems.

\end{abstract}

\maketitle

\setcounter{tocdepth}{1}
\tableofcontents

\newsect[INTR]{Introduction}

{By an {\slsf analytic cover} we understand a triple $(\tilde{X},\c,X)$ where $\c:\tilde X \rightarrow X$ is a zero-dimensional proper surjective holomorphic map between normal complex spaces $\ti X$ and $X$.

\noindent
{\slsf Zero-dimensional} means that for every $x\in X$ the preimage $\c ^{-1}(x)$ is discrete.
Since $\c$ is in addition supposed to be proper it is a {\slsf finite} map.
There exists a proper analytic set $\calr \subset X$ such that the restriction of $\c$ over $X\setminus \calr$ induces a finite unramified cover, see details in section \ref{STEIN_def}.
The set $\calr$ is called the {\slsf ramification locus} of the cover.
One says that  $(\tilde{X},\c,X)$ is a {\slsf connected analytic cover} if the covering space $\ti X$ is connected. 

The goal of this paper is to study when an analytic cover over a domain $D_0 \subset \cc^n$ can be extended to an analytic cover over a bigger domain $D_1$.
Clearly in order to achieve this it requires that the ramification locus extends to an analytic set in $D_1$.
The main problem is to understand if it is sufficient.
Examples are an important part of this paper.}
\newprg[INTR-results]{Statement of results}

\begin{nnthm}
\label{main_thm} 
Let $(\tilde D_0, \c_0 , D_0)$ be a $b_0$-sheeted {connected} analytic cover over a domain $D_0 \subset \cc^n$. Suppose that its ramification locus $\calr_0$ extends to an analytic set $\calr_1$ in a domain $D_1\supset D_0$ and the natural homomorphism $\isl_* : \pi_1(D_0\backslash\calr_0,z_*)\to\pi_1(D_1\backslash\calr_1,z_*)$ is surjective.
Here $z_*$ is some point in $D_0\backslash\calr_0$ and $\isl :D_0\to D_1$ denotes the natural inclusion. 
Then $(\tilde D_0, \c_0 , D_0)$ extends to a {connected} analytic
cover $(\tilde D_1, \c_1 , D_1)$ over $D_1$ in the  sense that there exists
a holomorphic map $\tilde{\isl}:\tilde D_0 \to \tilde D_1$ such that the diagram

\begin{equation}\label{diag_D}
\xymatrix{
    \tilde D_0 \ar[r]^{\tilde\isl} \ar[d]^{\c_0} & \tilde D_1 \ar[d]^{\c_1} \\
    D_0 \ar[r]^{\isl} & D_1
  }
\end{equation}
is commutative \ie $\c _1 \circ \tilde{\isl}=\c _0$ .
Moreover the cover $(\tilde D_1,\c _1,D_1)$ verifies the following:
\begin{enumerate}[label={\emph{(\roman*)}}]
\item \label{b_1<b_0 }the number $b_1$ of its sheets doesn't exceed $b_0$ and $\tilde{\isl}$ is surjective over $D_0$ \ie $\tilde{\isl}(\tilde D_0)=\c^{-1}_1(D_0)$;

\item \label{b_1=b_0} if $b_1=b_0$ then the map $\tilde{\isl}:\tilde D_0 \to \tilde D_1$ is injective;

\item \label{b'_1<b_1} every {connected} analytic cover $(\tilde D'_1,\c' _1,D_1)$ which is an extension of $(\tilde D_0,\c _0,D_0)$ has the number $b'_1$ of sheets not more than $b_1$;

\item \label{b'_1=b_1} if $b'_1=b_1$ then $(\tilde D'_1,\c' _1,D_1)$ and $(\tilde D_1,\c _1,D_1)$ are equivalent.

\end{enumerate}

\end{nnthm}

\smallskip
Two analytic covers $(\tilde X_0,\c_0 ,X)$ and $(\tilde X_1,\c_1,X)$ over a normal complex space $X$ are {\slsf equivalent} if there exists a biholomorphic map $\Phi:\tilde X_0 \to \tilde X_1$ such that $\c_1\circ \Phi = \c_0$.

The map $\tilde{\isl}$ of \refthm{main_thm} is surjective over $D_0$ 
but
may not necessarily be an injection. 
In general let $(\tilde X_0,\c_0,X_0)$ be an analytic cover over a normal complex space $X_0$. 
One says that it extends to an analytic cover $(\tilde X_1, \c_1,X_1)$ over a normal complex space $X_1$ if there exists a holomorphic imbedding $\isl :X_0\to X_1$ and a holomorphic map/injection $\tilde\isl :\tilde X_0\to \tilde X_1$ such that the
diagram 
\begin{equation} \label{diag_X}
\xymatrix{
    \tilde X_0 \ar[r]^{\tilde\isl} \ar[d]^{\c_0} & \tilde X_1 \ar[d]^{\c_1} \\
    X_0 \ar[r]^{\isl} & X_1
  }
\end{equation}
\noindent is commutative. 
Depending on whether $\tilde\isl$ is supposed to be an
injection or not one gets different notions of extension of analytic
covers. 
\refthm{main_thm} provides us an extension in a  ``weak" sense.
The conclusion of \refthm{main_thm} is not necessarily true if the morphism
\[ \isl : \pi_1(D_0 \backslash \calr_0,z_*) \to \pi_1(D_1 \backslash \calr_1,z_*)\]
is not surjective as it is shown by the following in \refsection{EXPLE_contre}.

\begin{nnexmp} \rm \label{exmp_braid}
There exist domains $D_0 \deff \Delta^4 \subset D_1=\cc^4$, a divisor $\calr_1 \subset D_1$ and a $3$-sheeted connected analytic cover $(\tilde D'_0,\c _0,D_0 \setminus \calr_1)$ such that
\begin{enumerate}[label={(\roman{*})}]
\item it does not extend to a $3$-sheeted connected analytic cover over $D_1 \backslash \calr_1$ ;

\item but it extends to a $4$-sheeted connected analytic cover $(\tilde D'_1,\c _1,D_1 \setminus \calr_1)$.
\end{enumerate}
\end{nnexmp}

\begin{nnrema}
In this paper speaking about extension of analytic covers we always assume that $\tilde{\isl}$ should be surjective over $X_0$, \ie $\tilde{\isl}(\tilde X_0)=\c_1^{-1}(X_0)$.
\end{nnrema}

\smallskip
For $R>0$ denote by $\Delta^k_R$ the $k$-dimensional polydisk of radius $R$ and set $\Delta^k_1=\Delta^k$.
The one dimensional disk of radius $R$ will be simply denoted $\Delta_ R$.
\begin{nnthm} \label{ext_disk}
Let $D_0\deff \Delta_R^{n-1}\times\Delta$ be a polydisk in $\cc^n$, $0<R<1$.
Let $(\tilde D_0, \c_0 , D_0)$ be a $b_0$-sheeted {connected} analytic cover over $D_0$. 
Suppose that its ramification locus $\calr_0$ extends to an analytic set 
$\calr_1$ in $\Delta^{n}$ such that $\calr_1$ does not intersect $\Delta^{n-1}\times\d\Delta$.
Then :
\begin{enumerate}[label=\emph{(\roman*)},ref={\thennthm}.{\emph{(\roman*)}}]
\item \label{ext_disk_b} $(\tilde D_0, \c_0 , D_0)$ extends in the weak sense to a $b_1$-sheeted {connected} analytic cover $(\tilde D_1, \c_1 ,\Delta^{n})$ over $\Delta^{n}$.
Moreover conditions \ref{b_1<b_0 }, \ref{b_1=b_0}, \ref{b'_1<b_1} and \ref{b'_1=b_1} in \refthm{main_thm} are satisfied.

\item \label{ext_disk_2} if $b_0=2$ then $(\tilde D_0,\c _0,D_0)$ uniquely extends in the strong sense to $(\tilde D_1,\c_ 1,\Delta^{n})$;  \ie the map $\tilde{\isl}$ in diagram \eqref{diag_D} is injective.

\end{enumerate}

\end{nnthm}

We prove also the Hartogs type extension theorem for analytic covers over a $2$-convex Hartogs figure.
Recall that a {\slsf $q$-convex Hartogs figure} is the following domain in $\cc^n$
\[
H^{n,n-q}_r=\Delta_r^{n-q}\times \Delta^q \cup \Delta^{n-q}\times\left(\Delta^q \backslash \overline{\Delta_{1-r}^q}\right)
\]
for some $0 < r <1$. 
\begin{nnthm}\label{ext_hart}
Let $(\tilde H,\c _0,H^{n,n-2}_r)$ be a $b_0$-sheeted {connected} analytic cover over a $2$-convex Hartogs figure ($n\geq 3$).
Then :
\begin{enumerate}[label=\emph{(\roman*)},ref={\thennthm}.{\emph{(\roman*)}}]
\item \label{ext_hart_b}the cover extends in the weak sense to a $b_1$-sheeted {connected} analytic cover over $\Delta^n$ which verifies conditions \ref{b_1<b_0 }, \ref{b_1=b_0}, \ref{b'_1<b_1} and \ref{b'_1=b_1} of \refthm{main_thm}.

\item \label{ext_hart_2} if $b_0=2$ then the cover $(\tilde H,\c _0,H^{n,n-2}_r)$ extends in the strong sense to $(\tilde D_1,\c_ 1,\Delta^n)$.

\end{enumerate}

\end{nnthm}

\begin{nnrema} \rm
\begin{enumerate}[label={(\alph{*})}]
\item There exists in \refthm{ext_disk_2} and \refthm{ext_hart_2} a holomorphic function $f\in\calo({\Delta^{n}})$ such that $\tilde D_1\deff \{(\zeta,z)\in\cc\times\Delta^{n} \tq \zeta^2=f(z)\}$ and $\c _1$ is induced by the canonical projection $\cc\times\Delta^{n} \to \Delta^{n}$.

\item The ``strong" parts of \refthm{ext_disk} and \refthm{ext_hart} are no longer true when the number of the sheets is larger than $2$ as it will be shown in the following example.
\end{enumerate}
\end{nnrema}

\begin{nnexmp} \rm \label{nnex_3} 
Set $\Omega \deff \left\{z_1 \in \cc  \tq \left|z_1-1\right|<1 \right\}$.
Let $f(z)$ be the branch of $\sqrt[3]{z}$ on $\Omega$ such that 
$f\left(1\right)=1$.
Let $(\tilde D_0,\c_0,D_0)$ be the $3$-sheeted {connected} analytic cover over $D_0\deff \Omega \times \cc$ defined as 
\[\tilde D_0 \deff \left\{ (z_,w)\in D_0\times\cc \tq w^3+\frac{f(z_1)}{\sqrt[3]{4}}w+\frac{iz_2}{\sqrt{27}}=0  \right\}\]
with $\c _0:\tilde D_0\to D_0$ induced by the natural projection.
Then there does not exist a {connected} extension of $(\tilde D_0,\c_0,D_0)$ over $D_1 \deff \cc^2$ with more than one sheet.
Remark that $(\tilde D_0,\c _0,D_0)$ can be extended to a $3$-sheeted analytic cover over the Hartogs figure 
$H\deff \left[\Omega\times\Delta_4\right]\cup \left[\Delta_4\times\left(\Delta_4\backslash\overline{\Delta_{3}}\right)\right]$.

\end{nnexmp}

{One says that an analytic cover $(\ti X,\c ,X)$ is {\slsf Galois} if the restriction over $X\setminus\calr$ induces a Galois regular cover.}

Let us give the following example of a non-extendible Galois analytic cover in the strong sense over the polydisk.
\begin{nnexmp}
\rm
\label{nnex_galois}
{In the notations of \refexmp{nnex_3} we set $D_0 = \Omega \times \cc$ and 
\begin{equation*}
\ti D_0 \deff \left\{  (z,w)\in D_0 \times \cc \tq w^3 - (z_1-z_2^2)(g({z_1})-z_2)=0   \right\},
\end{equation*}
where $g(z)$ is the branch of $\sqrt{z}$ on $\Omega$ such that $g(1)=1$.
Let $\c_0 : \ti D_0 \to D_0$ be the restriction of the canonical projection. 
Then $\ti D$ inherits a structure of a normal complex space such that $(\ti D_0,\c _0,D_0)$ becomes a $3$-sheeted Galois analytic cover over $D_0$ and there does not exist a connected analytic extension over $D_1 = \cc^2$ with more than one sheet.}
\end{nnexmp}

\newprg[INTR-struct]{Structure of the paper}

We begin by recalling in \refsection{STEIN} the definition of an analytic cover over a normal complex space and list some elementary properties of such covers.
Then we give a topological extension result which allows us to extend topologically a regular cover through an analytic hypersurface, see \refthm{thm_stein}.
This result is due to K. Stein, see \cite{St}.
A complete proof of this theorem is given in \refsection{STEIN_result} for the reader's convenience following the construction given by \cite{Ni}.
\refsection{GRA_REM} is about \refthm{thm_gra_rem} of H. Grauert and R. Remmert.
We recall in that section a complete proof of this Theorem using $L^2$ estimates methods given by G. Dethloff in \cite{De}, see \refthm{thm_L2}.
Then one obtains \refthm{struc_normal} which provides an analytic structure of the cover space.
One deduces in \refsection{THUL} some extension results for analytic covers, see \refthm{ext_norm_stein}, as well as one more example of an extension the so called Thullen-type extension theorem for analytic covers. 
We give it here for the sake of the future references.

\begin{nnthm} 
\label{ext_thullen} 
Let $X$ be a normal complex space, $A$ be a proper analytic subset of $X$ and $X_0$ be an open subset of $X$ which contains $X\setminus A$ and intersects every one-codimensional branch of $A$. 
Let $\c_0: \tilde X_0 \rightarrow X_0$ be an analytic cover over $X_0$. 
Then it uniquely extends in the strong sense to an analytic cover $\c:\tilde X\rightarrow X$ over the whole of $X$.
\end{nnthm}

We prove in \refsection{MAIN} the main result of this paper \ie \refthm{main_thm} and give some immediate consequences.
\refsection{POLYDISK} is about \refthm{ext_disk}.
\refthm{pic_sim} states that the hypothesis of \refthm{main_thm} are satisfied as soon as 
\begin{equation}\label{cs_polydisk}
\calr_1\cap \left(\Delta^{n-1}\times\d \Delta\right)=\emptyset.
\end{equation}
It is a somewhat precise version of \lemma{pic_sim_0} due to E. Picard and E. Simart.
We check that the proof of this Lemma given by T. Nishino in \cite{Ni} can be adapted to prove also \refthm{pic_sim}.
Then we prove \refthm{ext_disk}.
In \refsection{HART} we give the proof of \refthm{ext_hart} using some techniques of exhaustion by $(n-2)$-convex domains.
At the end of this section we give details about {\sl Examples} \ref{exmp_braid}, \ref{nnex_3} and \ref{nnex_galois}.

\smallskip 

I am grateful to  S. Ivashkovich and S. Orevkov who gave me valuable hints and examples for the proofs in this paper.

\newsect[STEIN]{Topological extension of covers}

\newprg[STEIN_def]{Analytic covers and their extensions}
Recall that a {\slsf regular cover} is a locally homeomorphic map
$\c:\tilde X\to X$ between Hausdorff topological spaces such that for
every $x_0\in X$ there exists a neighbourhood $U\ni x_0$ such that its preimage
$\c^{-1}(U)$ is at most countable disjoint union of its connected
components $\tilde{U}_i$ and for every $i$ the restriction $\restriction{\c}
{\tilde{U}_i}:\tilde{U}_i\to U$ is a homeomorphism. As it is well known (and 
obvious) if $\c :\tilde X\to X$ is a regular cover then for every path
$\gamma : [0,1]\to X$ and every $a\in \tilde X$ such that $\c (a) =
\gamma (0)$ there exists a unique lift $\tilde \gamma$ of $\gamma$ starting
at $a$, \ie a path $\tilde\gamma :[0,1]\to \tilde X$ such that
$\tilde\gamma (0)=a$ and $(\c\circ\tilde\gamma )(t) =\gamma (t)$ for
all $t\in [0,1]$.

\smallskip 

A regular cover $\c :\tilde X \to X$ is {\slsf finite} if it is proper.
In this case there exists $b\in\mathbb{N}$ such that every $x\in X$ has exactly $b$ preimages.
Recall the following facts, see \cite{Ha} for more details.
\begin{rema}\rm \label{corresp_cover}
Let $(\tilde X,\c ,X)$ be a $b$-sheeted regular cover.
Fix $z_*\in X$ and let $w_*\in \tilde X$ be a preimage of $z_*$ by $\c $.
We let $\pi_1(X,z_*)$ denote the fundamental group of $X$ with base point $z_*$.
For every path $\gamma \in \pi_1(X,z_*) $ we let $ \c ^* \gamma$ denote the lifted path of $\gamma$ starting at $w_*$.

\noindent{\bf(1).}\; Set $K \deff \{ [\gamma]\in \pi_1(X,z_*) \tq  \c ^* \gamma \mbox{ is closed in } \tilde X\}$.
It is a subgroup of index $b$ in $\pi_1(X,z_*)$.
Indeed let ${\alpha}$, ${\beta}$ be two closed paths starting at $z_*$ and let $\tilde{\alpha}\deff \c ^* \alpha$ and $\tilde\beta\deff \c ^*{\beta}$.
Then $\tilde{\alpha}(1)=\tilde{\beta}(1)\Longleftrightarrow [\beta\cdot\alpha^{-1} ]\in K$, which means that there are as many left cosets of $K$ in $\pi_1(X,z_*)$ as preimages of $z_*$ by $\c$.
Conversely if $K < \pi_1(X,z_*)$ is a subgroup of index $b$ there exists a finite regular cover $\c : \tilde X \to X$ with $b$ sheets and a fixed preimage $w_*$ of $z_*$ such that $K \deff \{ [\gamma]\in \pi_1(X,z_*) \tq  \c ^* \gamma \mbox{ is closed in } \tilde X \}$.
In other words $K=\im \c_*$ where $\c _*:\pi_1(\tilde X,w_*)\to \pi_1(X,z_*)$ denotes the natural monomorphism induced by $\c$.
\smallskip

\noindent{\bf(2).}\; Suppose that $\tilde{X}$ is path-connected and $X$ is path-connected, locally path-connected and locally simply connected. 
For every $\gamma \in \pi_1(X,z_*)$ and for every preimage $w_*$ of $z_*$ the ending point $\c^*(\gamma)(1)$ lies in the fiber of $z_*$.
Thus one obtains a homomorphism $\rho: \pi_1(X,z_*)\to {S}_b$ to the symmetric group of $b$ elements.
The map $\rho$ is called the {\slsf monodromy representation} of the cover.
Conversely suppose there exists a group homomorphism $\rho : \pi_1(X,z_*)\to {S}_b $ with {\slsf transitive image} \ie for every $i,j \in \{1,\dots,b\}$ there exists $\gamma \in \pi_1(X,z_*)$ which verifies $\rho(\gamma)(i)=j$.
Then there exists a $b$-sheeted regular cover $(\tilde X,\c _0,X)$ such that its monodromy representation coincides with $\rho$.
\end{rema}

\smallskip

A {\slsf regular analytic cover} of complex spaces is by definition a regular 
cover of normal complex spaces. It should be said that regular covers are
particular case of \rie{} domains. 

\begin{defi} 
\label{def_cover}
An {\slsf analytic cover} is a triple 
$(\tilde{X},\c,X)$ where $\c:\tilde X \rightarrow X$ is a zero-dimensional proper surjective holomorphic mapping between normal complex spaces $\ti X$ and $X$.
\end{defi}

\noindent
Let $\tilde{\calr} \subset \ti X$ be the set of critical points of $\c $.
It is a proper analytic set according to the Theorem of Sard on normal complex spaces, see for example \cite{Man}.
Since $\c $ is proper, Theorem of Remmert implies that its image $\calr \deff \c (\tilde{\calr})$ is proper analytic in $X$, see Satz 23 in \cite{Re}.
The restriction $\c : \tilde{X} \setminus \tilde{\calr} \to X\setminus \calr$ is locally biholomorphic.
By the hypothesis there exists $b\in \nn$ such that every point $x\in X\setminus \calr$ has exactly $b$ preimages.
It follows that $\c : \tilde{X} \setminus \tilde{\calr} \to X\setminus \calr$ induces a $b$-sheeted regular analytic cover.
Since $\ti \calr$ is normal the analytic set $\ti \calr$ does not locally separate it.
The set $\tilde\calr$ is called the {\slsf branching locus} of the cover and $\calr$ 
the {\slsf ramification locus}. 

\smallskip

{If $X\subset \cc^n$ is a domain we may assume $\calr$ is empty or a pure one codimensional analytic set since there cannot exist some branching point over any point of $\calr$ of an at least two codimensional component of $\calr$.

In general the ramification locus of an analytic cover over a normal complex space may have a codimension at least equal to $2$ as it is shown in the following example.

\begin{exmp}
Let $X\subset \cc^3$ be the analytic hypersurface defined by $z^2-xy=0$.
Since $\Sing(X)=\{0\}$ has codimension equal to $2$ in $X$ Oka's Theorem implies that $X$ is a normal complex space, see \cite{Be}.
Let $\c : \cc^2 \mapsto X$ be the holomorphic map given by $\c(s,t)=(s^2,t^2,st)$.
Then $(\cc^2,\c, X)$ is a $2$-sheeted analytic cover with ramification divisor $\calr = \{0\}$.
\end{exmp}
%

\begin{defi}
One says that an analytic cover  $(\ti X,\c ,X)$ is {\slsf connected} if $\ti X$ is connected.
\end{defi}
}

\newprg[STEIN_result]{Theorem of Stein}

The following statement can be found in \cite{St}, see Satz 1.
We give here the complete proof for the reader's convenience. 
We say that a closed subset $\calr$ of  a topological space $X$ doesn't locally separate it if for every $x_0\in \calr$ there exists a neighbourhood basis $\{U_{\alpha}\}$ of $x_0$ such that $U_{\alpha}\setminus \calr$ is connected.

\begin{thm} 
{}\label{thm_stein}
Let $X$ be a locally compact and locally connected Hausdorff topological space.
Suppose that $\calr$ is a proper closed subset of $X$ which does not 
locally separate it and let $\c_0:\tilde X_0 \rightarrow X\backslash \calr$ be a finite 
regular cover. 
Then there exists a unique locally compact, locally connected Hausdorff
topological space $\tilde X_1$, a unique imbedding $\tilde\isl :\tilde X_0 \rightarrow \tilde X_1$ surjective over $X\backslash\calr$
and a unique continuous surjective proper zero-dimensional map $\c_1 : \tilde X_1 \rightarrow X$ 
such that $\tilde \calr \deff\c_1^{-1}(\calr)$ is proper, closed, does not locally separate $\tilde X_1$ and the following diagram 

\begin{equation}  \label{ext_X_reg}
\xymatrix{
    \tilde X_0 \ar@{^{(}->}[r]^{\tilde\isl} \ar[d]^{\c_0} & \tilde X_1 \ar[d]^{\c_1} \\
    X\backslash\calr \ar@{^{(}->}[r]^-{\isl} & X
  }
\end{equation} 
is commutative \ie $\c_1 \circ \tilde{\isl}=\c _0$. 
\end{thm}
\proof
Without loss of generality we can suppose that $X$ and $\tilde X_0$ are connected. 
Let $b$ be the number of the sheets of $\c_0$.
For every $p\in \calr$ we use like in \cite{Ni} a basis $\{U_{\alpha}\}$ of connected neighbourhoods of $p$ such that $U_{\alpha}\backslash\calr$ is connected.
Every $\c_0^{-1}(U_{\alpha}\backslash\calr)$ can be decomposed into a finite number of connected components. 
If for every $\alpha$ there exists one of those components ${\delta}_{\alpha}$ such that
$${\delta}_{{\alpha}+1} \subset {\delta}_{\alpha}\mbox{ and } \bigcap_{\alpha} {\delta}_{\alpha} = \emptyset$$
we say that the sequence $\{{\delta}_{\alpha} \}_{\alpha}$ defines a {\slsf boundary point} $\tilde p$ of $\tilde X_0$ over $p$.
The sequence $\{\delta_{\alpha}\}_{\alpha}$ is a fundamental system of $\tilde p$ in $\tilde{X}_0$. 
Two systems $\left\{{\delta}_{\alpha} \right\}_{\alpha}$ and $\left\{{\eta}_{\beta} \right\}_{\beta}$ define the same boundary point
$\tilde p$ if for every ${\alpha}$ (resp. $\beta$ ) there exists some $\beta$ (resp. $\alpha$) such that ${\eta}_{\beta} \subset {\delta}_{\alpha}$ (resp. ${\delta}_{\alpha}\subset {\eta}_{\beta} $ ).
Set $q\in\calr\cap U_{\alpha} $ and let $\tilde q$ be a boundary point of $\tilde X_0$ over $q$. 
We say that $\tilde q$ {\slsf touches} $\delta_{\alpha}$ if there exist $\beta_0$ and a fundamental 
system $\left\{\eta_{\beta}\right\}_{\beta\geq \beta_0}$ of $\tilde q$ in $\tilde X_0$ which is contained in $\delta_{\alpha}$. 
Let $\tilde \calr$ be the set of boundary points over $\calr$ and set $\tilde X_1 \deff 
\tilde X_0 \cup\tilde\calr$.
For every ${\alpha}$ we let $\tilde\delta_{\alpha}$ denote the union of $\delta_{\alpha}$ with the set of points $\tilde q$  which touch $\delta_{\alpha}$ where $q\in\calr\cap U_{\alpha}$. 
The sequence $\left\{ \tilde\delta_{\alpha}
\right\}_{\alpha}$ is called a {\slsf fundamental neighbourhood system} of $\tilde p$ in $\tilde X_1$.
Remark that there are at most $b$ boundary points over $p$.  
 
\smallskip $\tilde X_1$ with the topology as above becomes a Hausdorff topological space.
Indeed let $\tilde{p_1}$ and $\tilde{p_2}$ be two different points in $\tilde X_1$.
Let us prove there exist two neighbourhoods of those points whose intersection is empty.
We shall prove this statement when $\tilde{p_1}$ and $\tilde{p_2}$ are different boundary points above the same point $p \in \calr$.
Let $\{{\delta}^1_{\alpha}\}_{\alpha}$ and $\{ {\delta}^2_{\alpha}\}_{\alpha}$ be their fundamental systems in $\tilde X_0$.
By the hypothesis we can suppose there exists some ${\alpha}_0$ such that for every ${\alpha}$ one has 
 ${\delta}^2_{\alpha} \not\subset {\delta}^1_{{\alpha}_0}$.
Since ${\delta}^1_{{\alpha}_0}$ and ${\delta}^2_{{\alpha}_0}$ are some connected components of $\c_0^{-1}(U_{{\alpha}_0})$ the fact that 
${\delta}^2_{{\alpha}_0} \not\subset {\delta}^1_{{\alpha}_0}$ implies that they are distinct and therefore disjoint. 
Moreover if $\tilde q \in \tilde\calr$ touches $\delta_{{\alpha}_0}^2$ 
then it can not touch $\delta_{{\alpha}_0}^1$. 
That's why $\tilde\delta_{{\alpha}_0}^1$ and $\tilde\delta_{{\alpha}_0}^2$ are two 
disjoint neighbourhoods of $\tilde p_1$ and $\tilde p_2$.
Moreover $\ti X_1$ is by construction locally connected. {Let us prove it is connected.
Let $E$ and $F$ be two disjoint open sets inside $\ti X_1$ such that $\ti X_1 = E \cup F$ and $\ti X_0 \subset E$. 
Suppose $F\neq \emptyset$.
Since $\ti X_1$ is locally connected one can find $\ti p \in \ti \calr \cap F$ and a connected neighbourhood $\ti U \supset \ti p$.
Then $\ti U \cap E$ and $\ti U \cap F$ are not disjoint and we have a contradiction.}

\smallskip

 We define $\c_1$ as follows. For every $\tilde p \in \tilde X_1$
\begin{equation}
  \label{ext_reg_cover}
\c_1(\tilde p)=
\left\{ \begin{array}{cc} 
\c_0(\tilde p) & \mbox{ if } \tilde p \in \tilde X_0\\
p & \mbox{ if } \tilde p\in\tilde\calr \mbox{ is a boundary point over }p.
\end{array}
\right. 
\end{equation}
The map $\c_1 : \tilde X_1 \rightarrow X$ is obviously continuous and surjective.
Let us prove it is open.
Let $\tilde U \subset \tilde X_1$ be an open set, $U \deff \c_1(\tilde U)$, $p \in U$ and $\tilde p \in \tilde U$ be a preimage of $p$ by $\c_1$.
If $p \notin \calr$ one can find an open set $\tilde V$ containing $\tilde p$ such that $\c_1$ is a homeomorphism on $\tilde V$.
Then $V \deff \c_1(\tilde V)$ is a neighbourhood of $p$ contained in $U$.
If $p \in \calr$ we take a neighbourhood basis $\{\tilde \delta_{\alpha}\}$ of $\tilde p$ such that $\tilde \delta_{\alpha} \subset \tilde U$.
Then $U_{\alpha}=\c_1(\tilde \delta_{\alpha}) \subset U$ is a neighbourhood of $p$ in $U$.
It follows that $U$ is open and we deduce the result.

Let us prove that $\c_1$ is closed.
Let $\tilde F$ be a closed subset of $\tilde X_1$ and take $p \in X \backslash F$ where $F\deff \c_1(\tilde F)$.
We firstly suppose that $p \notin \calr$.
Let $\tilde p \in \tilde X_1 \backslash \tilde F$ be a preimage of $p$ by $\c_0$.
By the hypothesis on $\tilde F$ one can find an open neighbourhood $\tilde U \subset \tilde X_1 $ of $\tilde p$ such that $\tilde U \cap \tilde F = \emptyset$ and also $U \cap F = \emptyset$ where $U \deff \c_1(\tilde U)$.
The fact that $\c_1$ is open implies that $U$ is an open subset of $X \backslash F$ containing $p$.
If $p \in \calr$ we take a neighbourhood basis $\{ U_{\alpha}\}$ of $p$.
Let us suppose that for every ${\alpha}$ we have $U_{\alpha} \cap F \neq \emptyset$.
Let $\tilde p$ be a preimage of $p$ by $\c_1$ and $\{\tilde \delta_{\alpha}\}_{\alpha}$ be its neighbourhood basis. 
Since $p\notin F$ one has $\tilde p \notin \tilde F$ and $\tilde \delta_{\alpha} \cap \tilde F$ is empty for some ${\alpha}$.
That gives us a contradiction because $U_{\alpha}\cap F \neq \emptyset$.
Therefore $\c_1$ is closed.

Moreover the fiber $\c_1^{-1}(p)$ of every $p\in X$ contains at most $b$ points by construction.
It follows that $\c_1$ is proper.
The set $\tilde \calr$ is clearly a proper closed subset of $\tilde X_1$ because $\calr$ is proper and $\c_1$ is continuous.
Moreover it does not locally separate $\tilde X_1$. 
Indeed let $\tilde p$ be a boundary point, $\tilde U$ be an open connected neighbourhood of $\tilde p$ and $\{\ti\delta_{\alpha}\}_{{\alpha}\geq {\alpha}_0}$ be a fundamental neighbourhood system of $\tilde p$ in $U$.
Since $\c_1(\tilde\calr)=\calr$ one obtains that $\tilde{\delta}_{\alpha} \backslash \tilde\calr=\delta_{\alpha}$ is connected.

Let us prove that $\tilde X_1$ is locally compact. 
Let $\tilde p \in \tilde X_1$ and let $\tilde K_{\alpha}$ be the connected component of $\c_1^{-1}(K_{\alpha})$ which contains $\tilde p$ where $K_{\alpha}\subset U_{\alpha}$ is a compact neighbourhood of $p$.
Then $\tilde K_{\alpha}$ is compact.
Indeed let $\tilde\calv\deff\{\tilde V_i \tq i\in I\}$ be an open covering of $\tilde K_{\alpha}$.
Since $\c_1$ is open $\calv\deff\{\c_1(\tilde V_i) \tq i\in I\}$ is a open covering of the compact set $K_{\alpha}$.
It follows there exists a finite subset $\{i_1,\dots,i_s\}$ of $I$ such that $\{\c_1(\tilde V_{i_1}),\dots,\c_1(\tilde V_{i_s})\}$ is an open covering of $K_{\alpha}$.
Then $\{\tilde V_{i_1},\dots,\tilde V_{i_s}\}$ is a finite open covering of $\tilde K_{\alpha} $ extracted from $\tilde\calv$.

\smallskip Finally let us prove the uniqueness of $\tilde X_1$.
 Let $\tilde{X}_1'$ be a locally compact Hausdorff topological space and $\c_1':\tilde{X}_1' \rightarrow X$ be a continuous proper zero-dimensional surjective map such that $\tilde{\calr}'\deff\c_1'^{-1}(\calr)$ is proper closed and does not locally separate $\tilde{X}_1'$. 
Suppose that there exist 
two imbeddings $\tilde\isl':\tilde X_0 \rightarrow \tilde X_1'$ and $\tilde\isl:\tilde X_0 \rightarrow 
\tilde X_1$ such that the following diagram
 
\begin{equation}
\label{2_ext_reg_cover}
\xymatrix{
\tilde X_1' \ar[d]_{\c'_1} & \tilde X_0\ar@{^{(}->}[r]^{\tilde\isl} \ar@{_{(}->}[l]_{\tilde\isl'} 
\ar[d]^{\c_0} & \tilde X_1 \ar[d]^{\c_1} \\
X & X\backslash\calr \ar@{^{(}->}[r]_{\isl} \ar@{_{(}->}[l]^{\isl}& X
  }
\end{equation}
\noindent is commutative \ie $\c_1\circ \tilde\isl=\c_0 =\c_1'\circ\tilde\isl'$.  We have to prove that 
there exists a homeomorphism $\Phi:\tilde X_1 \rightarrow \tilde X_1'$ such that $\Phi \circ \tilde\isl=
\tilde\isl'$.
 
\smallskip Let $\tilde z\in\tilde\isl(\tilde X_0)$ and $\zeta$ be its unique preimage in $\tilde X_0$ by 
$\tilde\isl$. 
The application $\Phi:\tilde\isl(\tilde X_0) \rightarrow \tilde\isl'(\tilde X_0)$ defined by $\Phi(\tilde z)=\tilde\isl'(\zeta)$ is homeomorphic and verifies $\Phi \circ \tilde\isl=\tilde\isl'$. 
Let us extend it continuously to $\tilde X_1$. 
Let $\tilde p\in \tilde\calr$, $\{\tilde\delta_{\alpha}\}_{\alpha}$ be a neighbourhood basis of $\tilde p$ such that $\tilde{\delta_{\alpha}}\backslash\tilde\calr$ is connected and $\tilde p_{\alpha} \in \tilde{\delta_{\alpha}}\backslash\tilde\calr$.
We denote $\tilde p'_{\alpha} \deff \Phi(\tilde p_{\alpha})$.
The fact that $(\tilde p_{\alpha})_{\alpha}$ converges to $\tilde p$ implies by continuity of $\Phi:\tilde X_1\backslash\tilde \calr \to \tilde X_1'\backslash\tilde \calr'$ that $(\c'_1\circ\Phi(\tilde p_{\alpha}) )_{\alpha}$ converges to $p\deff \c_1(\tilde p)$.
It follows by connectedness of $\tilde{\delta_{\alpha}}\backslash\tilde \calr$ that there exists a unique preimage $\tilde p'$ of $p$ by $\c'_1$ such that $\lim_{\alpha} \tilde p'_{\alpha}=\tilde p'\in \calr'$.
Such $\ti p'$ does not depend of the choice of $\tilde p_{\alpha}$.
We define $\Phi(\tilde p)\deff \tilde p'$.
One obtains a well-defined map $\Phi:\tilde X_1 \to \tilde X_1'$ which is continuous and bijective by construction.
\refthm{thm_stein} is proved.

\smallskip\qed

\begin{rema} \rm 

\begin{enumerate}[label=(\alph*)]
\item Taking as $X$ a domain in $\cc^n$ and $\calr$ a divisor in $X$ we see that \refthm{thm_stein} provides a {\slsf topological} extension of a finite regular cover across an analytic set.

\item \label{stein_homeo} Remark that if in the conditions of \refthm{thm_stein} $\c_0:\tilde X_0 \rightarrow X\backslash \calr$ 
is a one-sheeted regular cover then the extended map $\c_1:\tilde X_1 \rightarrow X$ is a homeomorphism. 

\item \label{stein_connected} {In the assumption of \refthm{thm_stein} the space $\ti X_1$ is connected if $\ti X_0$ is connected.}

\end{enumerate}
\end{rema}

\newsect[GRA_REM]{Theorem of Grauert and Remmert}

In this section the space $\cc^n$ is equipped with the norm $\|z\|=\max\{ |z_j|\tq j=1,\dots,n\}$.
Let $D$ be a bounded pseudoconvex domain in $\cc^n$, $\calr$ be a pure one-codimensional analytic subset of $D$, $\tilde D$ a locally compact Hausdorff connected topological space and $\c : \tilde D \rightarrow D$ a continuous proper mapping.
Set $Y\deff D\backslash \calr$ and $\tilde Y \deff \c ^{-1}(Y)$.
Suppose that $\tilde \calr \deff \c^{-1}(\calr)$ doesn't locally separate $\tilde D$ and $\restriction{\c}{\tilde Y}:\tilde Y \to Y$ is a finite regular cover.
Let $b$ denote the number of its sheets. 
$\tilde Y$ inherits a canonical holomorphic structure which makes $\restriction{\c}{\tilde Y}:\tilde Y \to Y$ holomorphic, see \cite{GR4}. 
By a {\slsf weakly holomorphic} function on $\tilde D$ one understands a 
continuous function which is holomorphic at univalent points of the cover, \ie at points of $\tilde Y$.

\newprg[GRA_REM_L2]{$L^2$ existence Theorem}

The following result is a $L^2$ estimate theorem for \rie{} domains which in the case of domains in $\cc^n$ is due to H\"{o}rmander, see Theorem 4.4.2 of \cite{Ho}.
For the reader's convenience we sketch here the proof given in \cite{De}.
Remark that a similar result is stated in \cite{NS}.

\smallskip

We set $\mathrm{d}\nu\deff\c^{\star}\left(\text{d}\lambda\right)$ where  $\text{d}\lambda\deff\left(\frac{i}{2}\right)^ n$d$z\,\wedge\,$d$\bar{z}$ is the Lebesgue measure on $\cc^n$.
Let $ \calc^{\infty}(\tilde Y)_{(0,1)}$ be the set of $(0,1)$-forms on $\tilde Y$ with $\calc^{\infty}$ coefficients, \ie
$g=\sum_{i=1}^{n}g_i d\bar{z_i}$ where $g_i\in \calc^{\infty}(\tilde Y)$.
$\calc_0^{\infty}(\tilde Y)_{(0,1)}$ denotes the subset of $(0,1)$-forms in $\calc^{\infty}(\tilde Y)_{(0,1)}$ which are compactly supported.
Set $|g|^2 \deff \sum_{i=1}^{n}|g_i|^2$ and set $\phi_0(z)= 2n\ln\|z\|$, $z\in\cc^n$.

\begin{thm}\label{thm_L2}
With the previous assumptions let $g \in \calc^{\infty}(\tilde Y)_{(0,1)}$ be such that
$\dbar g=0$ and $\lambda_0\deff\int_{\tilde Y}|g|^2 e^{-\phi_0\circ\c} \mathrm{d}\nu <+\infty$.
Then there exists $u\in \calc^{\infty}(\tilde Y)$ such that $\dbar u=g$ and $ \int_{\tilde Y}|u|^2 e^{-\phi_0\circ\c} \mathrm{d}\nu \leq k \lambda_0$ where $k\deff(1+(\diam Y)^2)^2$.

\end{thm}

\proof
Let $\phi:Y\to \rr\cup\{-\infty\}$ be an upper semi-continuous function on $Y$.
We define the following vector spaces.

\[L^2(Y,\phi)=\left\{f:Y \to \cc \tq \int_Y |f|^2 e^{-\phi}  d\lambda < + \infty \right\} \]
and
\[ L^2(\tilde Y,\phi\circ\c )=\left\{f:\tilde Y \to \cc \tq \int_{\tilde Y} |f|^2 e^{-\phi\circ \c}  \mathrm{d}\nu < + \infty \right\}.\]
In the same way we define $L^2(Y,\phi)_{(p,q)}$ and $L^2(\tilde Y,\phi\circ\c )_{(p,q)}$ the set of $(p,q)$-forms whose coefficients belong to $L^2(Y,\phi)$ and $L^2(\tilde Y,\phi\circ\c)$, respectively.
Let $\phi_i$, $i=1,2,3$ be continuous functions.
An element $u \in L^2(\tilde Y,\phi_1\circ\c)$ is in $D_T$ if and only if $\dbar u \in L^2(\tilde Y,\phi_2\circ \c)_{(0,1)}$.
Since $D_T$ contains $\calc_0^{\infty}(\tilde Y)$ which is dense in $L^2(\tilde Y,\phi\circ\c )$ one can densely define the operator $T=\dbar:L^2(\tilde Y,\phi_1\circ\c)\to L^2(\tilde Y,\phi_2\circ\c)_{(0,1)}$ on Hilbert spaces.
In the same way one can define $D_S$ and the operator $S=\dbar:L^2(\tilde Y,\phi_2\circ\c)_{(0,1)}\to L^2(\tilde Y,\phi_3\circ\c)_{(0,2)}$.
We let $T^*$ denote the adjoint operator of $T$.
Let $\phi \in \calc^{\infty}(Y)$ be strictly plurisubharmonic and $\gamma:Y\mapsto\rr^*_+$ be continuous such that
\begin{equation} \label{L2_gamma}
\forall w\in \cc^n \text{ one has } H_z(\phi,w)\geq \gamma(z)\|w\|^2,
\end{equation}
where $H_z(\phi,w)$ denotes the Levi form of $\phi$.
Let $\{K_{\nu}\}_{\nu}$ be an exhaustion of $Y$ by compact sets, $\eta_{\nu}$ be a cutoff function in $\calc_0^{\infty}(Y)$ such that $0 \leq \eta_{\nu}\leq 1$ and $\restriction{\eta_{\nu}}{K_{\nu}}\equiv 1$.
Let $\psi \in \calc^{\infty}(Y)$ be a smooth function such that $\sum_{k=1}^{n}|\partial\eta_{\nu}/\partial\bar{z}_k|^2 \leq e^{\psi}$ on $Y$ and $\phi_i\deff \phi + (i-3)\psi$ $(i=1,2,3)$.

The first result needed to prove \refthm{thm_L2} is the following inequality. For every $ f \in \calc_0^{\infty}(\tilde Y)_{(0,1)}$ one has

\begin{equation}\label{L2_long}
 \int_{\tilde Y}\left( \gamma\circ\c -2 |\partial(\psi\circ\c)|^2 \right)|f|^2 e^{-\phi\circ\c }\mathrm{d}\nu \leq 2 \|T^{*}f\|^2_{\phi_1\circ\c}+\|Sf \|^2_{\phi_3\circ\c},
\end{equation}
see \cite{De} Lemma 2.1.
We let $A$ denote the subset of $(0,1)$ forms in $D_{T^*} \cap D_S$ which are compactly supported.
One deduce the following main result.
\begin{lem}
$\calc^{\infty}_0(\tilde Y)_{(0,1)}$ is dense  in $D_{T^*} \cap D_S$ for the graph norm
\[f \mapsto \|f \|_{\phi_2\circ\c}+ \|T^*f \|_{\phi_1\circ\c} + \|Sf \|_{\phi_3\circ\c}.\]
\end{lem}
The proof consists to see that $A$ is dense in $D_{T^*} \cap D_S$ and $\calc^{\infty}_0(\tilde Y)_{(0,1)}$ is dense in $A$ for the graph norm.
One mainly uses Hahn-Banach Theorem, Riez representation Theorem and Lebesgue's dominated convergence Theorem.
By using the existence of a smooth strictly plurisubharmonic exhaustion function $s$ on the pseudoconvex domain $Y$ one obtains the following result, see Lemma 4.4.1 of \cite{Ho} for more details.
\begin{lem}\label{L2_bounded}
Suppose that the function $\gamma$ defined in $\eqref{L2_gamma}$ is bounded and its lower bound is strictly positive. Then for every $g \in L^2(\tilde Y,\phi\circ \c)_{(0,1)}$ such that $\dbar g=0$ there exists $u \in L^2(\tilde Y,\phi\circ \c)$ which verifies $\dbar u=g$ and
\[\int_{\tilde Y} |u|^2 e^{-\phi\circ\c} \mathrm{d}\nu \leq 2 \int_{\tilde Y} \frac{|g|^2}{\gamma\circ\c} e^{-\phi\circ\c} \mathrm{d}\nu < +\infty.  \]
\end{lem}
We assume now that $\phi$ is a $\calc^{\infty}$ plurisubharmonic function on $Y$. 
We can apply the previous Lemma to {$z\mapsto \phi(z) + 2\ln (1+\|z\|^2)$} which is strictly plurisubharmonic 
and \mbox{$\gamma(z)\deff (1+\|z\|^2)^{-2}$}. 
One obtains the following statement.
\begin{lem} \label{L2_relax}
There exists $u \in L^2(\tilde Y,\phi\circ\c)$ such that $\dbar u = g$ and 
\[\int_{\tilde Y} |u|^2 e^{-\phi\circ\c}(1+\|\c\|^2)^{-2} \mathrm{d}\nu \leq \int_{\tilde Y}|g|^2 e^{-\phi\circ\c} \mathrm{d}\nu.\]
\end{lem}

Now let us prove \refthm{thm_L2}.
Let $a\in\rr$.
Set $Y_a\deff \{s(z)<a\}$ and $\tilde Y_a \deff \c^{-1}(Y_a)$.
According to Theorem 2.6.3 of \cite{Ho} there exists $\phi_a$ a $\calc^{\infty}$ plurisubharmonic function on $Y_a$ such that $\phi_{a}\underset{a \to +\infty}{\searrow} \phi_0$.
By \lemma{L2_relax} there exists $u_a \in L^2(\tilde Y_a,\phi\circ\c)$ such that
\[\int_{\tilde Y_a} |u_a|^2e^{-\phi_a \circ\c}(1+\|\c\|^{2})^{-2}\mathrm{d}\nu \leq 
\int_{\tilde Y_{a}}|g|^2 e^{-\phi_a\circ\c} \mathrm{d}\nu \leq \lambda_0.
\]
Then there exists a subsequence $(u_{a_j})_j$ from $(u_a)$ which weakly converges on every $\tilde Y_{a_j}$ to a function $u\in L^2(\tilde Y,\phi_0\circ\c)$ such that
\[\forall j\in \nn^* \text{ one has }\int_{\tilde Y_{a_j}} |u|^2e^{-\phi_{a_j} \circ\c}\left(1+\|\c\|^{2}\right)^{-2}\mathrm{d}\nu \leq \lambda_0.
\]
Set $k\deff(1+(\diam Y)^2)^2$. 
By a Cauchy-Schwarz inequality one obtains
\[\begin{array}{ccc}
\left(\int_{\tilde Y_{a_j}}|u|^2 e^{-\phi_{a_j}\circ \c}\textrm{d}\nu \right)^2 &\leq&
\left(\int_{\tilde Y_{a_j}}\frac{|u|^2}{(1+\|\c\|^2)^2}e^{-\phi_{a_j}\circ \c}\textrm{d}\nu \right)\left( \int_{\tilde Y_{a_j}}|u|^2(1+\|\c\|^2)^2e^{-\phi_{a_j}\circ \c}\textrm{d}\nu \right) \\
&&\\
&\leq& \lambda_0 \, k \, \left(\int_{\tilde Y_{a_j}}|u|^2e^{-\phi_{a_j}\circ \c}\textrm{d}\nu \right)
\end{array} 
 \] 
and then
\[ \int_{\tilde Y_{a_j}}|u|^2 e^{-\phi_{a_j}\circ \c}\textrm{d}\nu \leq k \, \lambda_0.
\]
One finally obtains by monotone convergence Theorem
\[ \int_{\tilde Y_{a}}|u|^2 e^{-\phi_{0}\circ \c}\textrm{d}\nu \leq k \, \lambda_0.
\] 
By the classical regularity of the $\dbar$-equation, see Theorem 4.2.5 of \cite{Ho} the solution $u$ is $\calc^\infty$ smooth and \refthm{thm_L2} is proved.

\smallskip \qed

\newprg[GRA_REM_thm]{Grauert-Remmert Theorem}

 The following theorem is due to Grauert-Remmert, see \cite{GR2} and \cite{Ni}. 
We use the approach of \cite{NS} and \cite{De} for the proof of Grauert-Remmert theorem concerning 
the existence of weakly holomorphic functions.

\begin{thm} {\slsf (Grauert-Remmert)}
\label{thm_gra_rem} Let $D\subset \cc^n$, $Y = D\backslash \calr$ and $\c : \tilde D \to D$ be as in \refthm{thm_L2}.
Then for every point $z_*\in Y$ with 
$\c^{-1}(z_*) = \{w_1,...,w_b\}$ there exists a weakly holomorphic function $h$ on
$\tilde D$ which takes at $\{w_1,...,w_b\}$ pairwise different
values.
\end{thm}
\proof Fix some $z_* \in D\backslash \calr$ and let $\c^{-1}(z_*)=\{w_1,\dots,w_b\}$ be as in the formulation of the theorem.
Since $\c:\tilde Y \to Y$ is locally biholomorphic one can find a function $p\in\calc^{\infty}_0(\tilde Y)$ such that $p\equiv 1$ in a neighbourhood of $w_1$ and $p\equiv 0$ in neighbourhoods of $w_2,\dots,w_b$.
One can deduce from \refthm{thm_L2} there exists $u\in\calc^{\infty}(\tilde Y)$ such that $\dbar u=\dbar p$, $u(w_i)=0$ $(i=1,\dots,b)$ and $\int_{\tilde Y} |u|^2 \mathrm{d}\nu < +\infty$.
Then $h_0\deff p-u$ is a holomorphic function on $\tilde Y$ such that
 $$h_0(w_i)=\left\{ \begin{array}{cc} 
                   1 & \mbox{ if }i=1 \\
                   0 & \mbox{ if }i\neq1,
                  \end{array}
 \right. $$
 and $\int_{\tilde Y}|h_0|^2 \mathrm{d}\nu < \infty$.
The fact that $D$ is a domain of holomorphy implies that there exists a holomorphic function $f$ on $D$ such that $f(z_*)\neq 0$ and $\restriction{f}{{\calr}}\equiv 0$, see \cite{GF}. 
In order to extend $h_0$ through the points of $\c^{-1}\left( \Reg(\calr) \right)$ one defines as in \cite{De} 
 $$t(w)=\frac{1}{f(z_*)}(f\circ \c)(w).$$ 
Then $t$ is a weakly holomorphic function on $\tilde D$ such that $\restriction{t}{\tilde \calr}\equiv 0$. 
Moreover for every $ i\in\{1,\dots,b\}$ one has $t(w_i)=1$ therefore $(t \cdotp h_0)(w_i)=h_0(w_i)$. 

Set $\tilde{Y}_0\deff\c^{-1}\left( D \backslash \Sing \calr \right) $. 
By Lemma 1.8 of \cite{De} $t\cdotp h_0$ extends to a weakly holomorphic function $h'$ on $\tilde{Y}_0$.

Now let us prove that $h'$ is locally bounded on $\tilde D$. 
We define the Weierstrass polynomial of $h'$:

\begin{equation}
 \label{pol_weier}
\forall \zeta\in \cc, \, \forall z\in D\backslash {\calr}, \; \omega_{h'}(z,\zeta)\deff \prod_{w\in \c^{-1}(z)}(\zeta-h'(w))= \zeta^ b + \sum_{i=1}^b a_i(z) \zeta^{b-i}.
\end{equation}
For every $i \in \{1,\dots, b\}$ the function $a_i$ is holomorphic on $D\backslash {\calr}$ and continuity of $h'$ on $\tilde Y_0$ implies that $a_i$ is continuous on $D \backslash \Sing({\calr})$. 
According to \rie{} Extension Theorem the function $a_i$ holomorphically extends to $D \backslash \Sing({\calr})$.
Since $\codim\left(\Sing({\calr})\right) \geq 2$ there exists a holomorphic function on $D$ which extends $a_i$. 
We still denote by ${a_i}$ this function.
One obviously has the following.
\begin{lem}
\label{localisation}
Let $P = \zeta^ b + \sum_{i=1}^ b a_i \zeta^ i$ be a complex polynomial. Then every root $\zeta_0$ of $P$ verifies the following 
inequality
 
\[
| \zeta_0 | \leq \max \left\{1 , \sum_{i=1}^ b |a_i| \right\}.
\]
 \end{lem}
\noindent According to the previous statement one has $|h'(w)| \leq \max\left\{1 , \sum_{i=1}^ b |a_i(\c(w))| \right\}$ for every $w \in \tilde{Y_0}$. The function $h'$ is also bounded near points of $\tilde{D}\backslash \tilde{Y_0}$ and therefore $h\deff t\cdotp h'$ is weakly holomorphic on $\tilde D$ and separates the sheets.
\refthm{thm_gra_rem} is proved.

\qed

\smallskip

Theorem of Grauert-Remmert provides an analytic structure on the {extended to\-po\-lo\-gi\-cal space}.
The main idea of the proof follows \cite{DG}. 

\begin{thm} \label{struc_normal}
In the conditions of \refthm{thm_gra_rem} the space $\tilde D$ inherits a unique structure of a normal complex space such that $\c$ becomes holomorphic, $\tilde\calr$ analytic and therefore $(\tilde D,\c,D)$ becomes an analytic cover. 
The structure sheaf of $\tilde{D}$ is the sheaf $\mathcal{O}'_{\tilde D}$ of weakly holomorphic functions.
\end{thm}
\proof 
 We shall prove that for every $w \in \tilde D$ we can find a neighbourhood $\tilde U$ of $w$ such that $\left(\tilde U,\mathcal{O}'_{\tilde U}\right)$ is normal.
Since $\tilde D\backslash\tilde{\calr}$ is locally homeomorphic to a domain of $D\backslash\calr$ it is sufficient to prove the result when $w \in \tilde \calr$.
Let $U$ be an open polydisk centered at $z\deff \c (w)$ in $D$ such that $\tilde\calr$ does not separate the connected component $\tilde U$ of $\c^{-1}(U)$ containing $w$ \ie $\tilde U \backslash \tilde\calr$ is connected.
 We can apply \refthm{thm_gra_rem} of Grauert-Remmert to the restriction $\restriction{\c}{\tilde U}:\tilde U 
\rightarrow U$ which we still denote by $\c$. 
Let $z_*$ be a fixed point of $U\backslash\calr$, $h$ be a holomorphic function on $\tilde U$ which separates the preimages of $z_*$, $\omega_h(z,\zeta)$ its Weierstrass polynomial on $U$ as in \ref{pol_weier} and $D(z)$ its discriminant.
Set $\sigma\deff\{z\in U\tq D(z)=0  \}$, $\tilde\sigma\deff\c^{-1}(\sigma)$ and $M\deff\left\{(z,\zeta)\in 
U\times \cc \tq \omega_h(z,\zeta)=0 \right\}$.
Consider the mapping
\[ \fonction{\Phi}{\tilde U}{U\times \cc}{x}{(\c(x),h(x))}. \]
Then $\Phi(\tilde U)=M$ is a connected analytic subset of $U\times \cc$.
Moreover $z\in \sigma$ if and only if $h$ does not separate some preimages of $z$ by $\c$. 
That's why the restriction $\Phi_1\deff\restriction{\Phi}{\tilde U \backslash \tilde\sigma}:\tilde U \backslash \tilde\sigma \rightarrow (U\backslash\sigma) \times \cc$ is injective. 
Mappings $\Phi$ and $\Phi_1$ are continuous, proper and surjective onto $M$ and $M\backslash\left(\sigma\times\cc \right)$ respectively. 
By Oka's normalization Theorem there exists a normal complex space $N$, a holomorphic map $\Psi:N \to M$ and an analytic set $A\subset \Sing M$ of $M$ such that
\begin{itemize}
\item $\Psi^ {-1}(A)$ is nowhere dense in $N$, 
\item $\Psi:\Psi^{-1}(A)\to A$ is finite and
\item $\Psi:N\backslash\Psi^{-1}(A)\to M\backslash A$ is biholomorphic.
\end{itemize}
Set $M_{\sigma}\deff M \cap (\sigma \times \cc )$ and consider the biholomorphism 
\[
t\deff\Phi^{-1}\circ \Psi: N\backslash \Psi^{-1}(\Ms)\to\tilde U \backslash\Phi^{-1}(M_{\sigma}).
\] 
The fact that $N$ is normal implies that $t:N\backslash \Psi^{-1}(\Ms) \to \tilde U \backslash\Phi^{-1}(M_{\sigma}) $ holomorphically extends to $\overline{t}: N \to \tilde U$.
Indeed let $w\in \Psi^{-1}(\Ms)$ and $\tilde W$ be an open neighbourhood of $w$ in $N$ such that $t:\tilde W\backslash\Psi^{-1}(\Ms)\to \tilde V$ is holomorphic where $\tilde V$ is biholomorphic to a connected analytic subset of some bounded domain in $ \cc^m$.
One can replace $t$ by $m$ holomorphic and bounded functions $t_i:\tilde W\backslash\Psi^{-1}(\Ms) \to \cc$ ($i=1,\dots,m$).
By normality of $\tilde W$ every $t_i$ extends to a holomorphic function $\overline{t_i}$ on $\tilde W$.
One obtains that the map $(\overline{t_1},\dots,\overline{t_m})$ gives an extension $\overline{t}:N \to \tilde U $ of $t$ on the whole of $N$.
Since $\Psi^{-1}(\Ms)$ and $\Phi^{-1}(\Ms)$ do not locally separate $N$ and $\tilde U$ respectively one can apply uniqueness of \refthm{thm_stein}  to the one-sheeted regular cover $(\tilde U \backslash\Phi^{-1}(M_{\sigma}) ,t, N\backslash \Psi^{-1}(\Ms))$.
It follows that $\overline{t}:N \to \tilde U$ is a biholomorphism and then $(\tilde U, \calo'_{\tilde U})$ and $(N,\calo_N)$ are isomorphic.
\smallskip\qed

\begin{exmp}\rm \label{nnex_stein}
Set $\tilde D \deff \{(z,w)\in\cc \times \cc \tq w^2=z \}$ and let $\c:\tilde D \to \cc $ be the map induced by the projection $(z,w) \mapsto z$.
Then $(\tilde D, \c, \cc)$ is a $2$-sheeted analytic cover with ramification divisor $\{z=0 \}$.
Let us define for every $(z,w)\in \cc^2$ the polynomial $h(z,w)=(z-1)w$.
Then $h \in \calo'(\tilde D)$ and separates the preimages of $z=-1$ since $h(-1,i)=-2i$ and $h(-1,-i)=2i$ but does not separate the preimages of $z=1$ because $h(1,1)=h(1,-1)=0$.
Let $\omega(z,\zeta)$ be the Weierstrass polynomial of $h$ defined in \ref{pol_weier} \ie 
$\omega(z,\zeta)=\left[ \zeta-(z-1)w \right]\left[ \zeta+(z-1)w \right]=\zeta^2-z(z-1)^2$ and
$M\deff \{(z,\zeta)\in \cc^2 \tq \omega(z,\zeta)=0 \}$.
Set $\sigma=\{z\in \cc \tq z(z-1)^2 = 0  \}=\{0,1
\}$.

One obtains that the map 
\[ \fonction{\Phi}{\tilde D}{M}{(z,w)}{(z,(z-1)w)}
\]
is well defined, continuous, surjective and the restriction of $\Phi$ over $M\backslash(\{1,0\})$ is injective.
Since $\Phi^ {-1}(1,0)=\{(1,1); (1,-1) \}$ the map $\Phi: \tilde D \to M $ is not globally injective. 
It follows that $(\tilde D, \Phi,M)$ and $(M,\id,M)$  are two different extensions of the cover $\left(\tilde D\backslash \{(1,1),(1,-1) \},\Phi, M\backslash\{(1,0)\}\right)$ over $M$. 
This happens because $\{(1,0)  \}$ locally separate $M$.
Indeed let $U$ be a small disk centered at $z=1$ which does not contain $0$ and $f(z)$ a branch of $\sqrt[2]{z}$ on $U$.
Then $M \cap (U \times \cc )$ can be decomposed into the union of the branches :
\[M_1\deff \{\zeta-f(z)(z-1)=0 \} \text{ and }M_2\deff\{\zeta+f(z)(z-1)=0 \}.
\]
which may intersect  when $f(z)(z-1)=0$ \ie at $(1,i)$.
Therefore $M \cap (U\backslash\sigma \times \cc )$ is not connected.
\end{exmp}

\newsect[THUL]{Thullen type extension of analytic covers}  

We give in this section some extension theorems for analytic covers. 
Let us start with the following result, see Proposition 3.3 in \cite{DG}.
\begin{lem}\label{GPR}
Let $\tilde X$ be a topological space, $X$ be a normal complex space and $\c:\tilde X \to X$ be a continuous map.
Then $(\tilde X ,\c , X)$ is an analytic cover if and only if for every $z\in X$ there exists 
an open neighbourhood $U \subset X$ of $z$ and 
an analytic cover $(U ,\pi ,V)$ over a domain $V \subset \cc^n$ such that 
the composed map $(\ti U,\pi\circ \c,V)$ is again an analytic cover 
where $\tilde U \deff \c^{-1}(U)$.
\end{lem}

\proof
The condition is obviously necessary. Let us prove the converse.
Let $\Sing X$ be the singular locus of $X$, $\tilde S \deff \c^{-1}(\Sing X))$ and $z \in X\backslash\Sing(X)$.
Let $(U , \pi , V)$ be an analytic cover as in the formulation of the lemma.
Since $(\ti U,\pi\circ \c,V)$ is an analytic cover the fact that $\c^{-1}(U)=(\pi \circ \c)^{-1}(V)$ implies that the restriction $ (\tilde X \backslash \tilde S ,\c, X \backslash \Sing X)$ is an analytic cover. 
We let $\calr$ denote its ramification locus.
Since $X$ is normal one has $\codim \Sing X \geq 2 $ and $\calr \subset X \backslash \Sing X$ uniquely extends to a proper analytic set $\overline{\calr}\subset X$ which does not locally separate it, see Theorem 9.4.2 in \cite{GR3}. 
It implies that $(\tilde X ,\c , X)$ is an analytic cover of ramification locus $\overline{\calr}$.

\qed

 \smallskip
 
A consequence of Grauert-Remmert and Stein Theorems is the following extension result for analytic covers.

\begin{thm}
\label{ext_norm_stein}
Let $(\tilde{X}_0 ,\c _0, X\backslash \calr)$ be a finite regular analytic cover where $X$ and $\tilde{X}_0$ are normal complex spaces and $\calr$ is a proper analytic subset in $X$ which does not locally separate it. 
Then there exists a unique normal complex space $\tilde X_1$, a unique imbedding $\tilde\isl:\tilde X_0 \rightarrow \tilde X_1$ and an analytic cover $  (\tilde X_1 ,\c _1, X)$ such that the diagram \eqref{ext_X_reg} is commutative \ie $\c_1\circ\tilde\isl=\c_0$. 
\end{thm}
\proof According to \refthm{thm_stein} there exist a unique locally compact Hausdorff topological space $\tilde X_1$,
a unique imbedding $\tilde{\isl}:\tilde X_0 \to \tilde X_1$ and a unique continuous surjective proper zero-dimensional map $\c_1:\tilde X_1 \rightarrow X$such that the diagram \eqref{ext_X_reg} is commutative and $\tilde \calr\deff\c^{-1}(\calr)$ does not locally separate $\tilde X_1$.
If $X$ is a domain of $\cc^n$, \refthm{struc_normal} implies that $\tilde X_1$ inherits a structure of normal complex space and the result follows.

We suppose now that $X$ is a normal complex space. 
Let $z\in X$.
There exist an open subset $U$ of $X$ containing $z$, an open polydisk $V\subset \cc^n$ and an analytic cover $(U ,\pi, V)$ of ramification divisor $\calr_V$.
According to \lemma{GPR} it is sufficient to prove that $(\tilde U ,\pi\circ\c_1, V)$ is an analytic cover where $\tilde U \deff \c_1^{-1}(U)$.
Let $\Sing(X)$ denote  the singular locus of $X$.
{Since $\pi$ is proper Theorem of Remmert implies that $\pi(\calr \cup \Sing(X))$ is analytic in $V$.
The set $\calr_1 \deff \calr_V \cup \pi(\calr \cup \Sing(X))$ is analytic in $V$ and does not locally separate it.
The map $\pi\circ\c_1$ induces also a finite regular cover over $V \setminus \calr_1$.}
According to the result proved above it uniquely extends to an analytic cover $(\tilde V, {\c}_V,V)$ over the whole of $V$. 
The fact that $\pi \circ \c_1:\tilde U \to V$ is a continuous surjective proper zero-dimensional map implies by uniqueness of \refthm{thm_stein} that there exists a homeomorphism $\Phi:\tilde U \to \tilde V$ such that $\pi \circ \c_1=\c _V\circ \Phi$.
That homeomorphism becomes a biholomorphism and it follows that $(\tilde U,\pi\circ\c_1 ,V)$ is an analytic cover.

\smallskip
\qed

\begin{rema} \label{grauert_connected} \rm
In the assumptions of \refthm{ext_norm_stein} the space $\ti X_1$ is connected if $\ti X_0$ is connected, see {\sl Remark 2.2 \ref{stein_connected}}.

\end{rema}
 
\begin{corol}
\label{Det_Gra}
Let $X$ be a normal complex space and $A\subset X$ be a proper analytic subset which does not locally separate it.
Let $(\tilde X_0 ,\c _0, X\backslash A)$ be an analytic cover with ramification locus $\calr \subset X\backslash A$.  
Assume that $\calr$ extends to an analytic subset $\overline{\calr}$ of $X$ which does not locally separate it.
Then $(\tilde X_0,\c _0, X\backslash A)$ uniquely extends to an analytic cover $(\tilde X_1,\c _1, X)$ over the whole of $X$.
\end{corol}

\proof 
We set $Y\deff X\backslash \overline{\calr}$ and $\tilde Y\deff\c_0^{-1}(Y)$.
The restriction $\restriction{\c_0}{\tilde Y}:\tilde Y \rightarrow Y$ induces a finite regular cover which uniquely extends to an analytic cover $(\tilde X_1,\c_1,X)$ according to \refthm{ext_norm_stein}.
One obtains two analytic covers $(\tilde X_0,\c _0,X\backslash A)$ and $(\c _1^{-1}(X\backslash A),\c _1,X\backslash A)$ which extend $(\tilde Y,\c _0, Y)$. 
By uniqueness they are equivalent and we can find a holomorphic imbedding $\Phi:\tilde X_0 \hookrightarrow \tilde X_1 $ such that $\c _1\circ\Phi=\c _0$. 
Then $(\tilde X_1,\c_1,X)$ also extends the cover $(\tilde X_0,\c_0,X\backslash A)$.

\smallskip\qed

\noindent{\slsf Proof of Thullen type extension.}
Now let us deduce the Thullen-type extension Theorem \ref{ext_thullen} from Introduction.
Let $\calr\subset X_0$ be the ramification locus of the cover $(\tilde X_0 ,\c_0, X_0)$.
According to the Thullen-Remmert-Stein Theorem (see \cite{Si}) the closure $\overline{\calr}$ of $\calr$ in $X$ is a proper analytic set which does not locally separate $X$.
Set $Y \deff X \backslash\overline\calr $ and $\tilde Y \deff \c_0^{-1}(Y)$.
Then $(\tilde Y,\c _0,Y)$ is a finite regular cover.
According to \refthm{ext_norm_stein} it uniquely extends to a cover $(\tilde X,\c ,X)$ over the whole of $X$.
By uniqueness of the extension it implies that $(\tilde X,\c ,X)$ extends $(\tilde X_0,\c_0,X_0)$.

\smallskip\qed

\begin{rema} \rm 

The fact that Thullen type extension holds for analytic covers 
follows from the fundamental papers of K. Stein \cite{St} and Grauert-Remmert \cite{GR2} is
known to the experts, see for example \cite{DG}. But a complete proof to our best knowledge
cannot be found in one place in the literature. 
Therefore we give in this paper a reasonably self-contained proof of the results stated above.
\end{rema}

\newsect[MAIN]{Extension by the extension of the ramification divisor}
In this section we prove Theorem \ref{main_thm} from the Introduction.
Set $Y_1\deff D_1\backslash\calr_1$, $Y_0=D_0\backslash \calr_0$, $\tilde{\calr}_0 \deff \c _0^{-1}(\calr_0)$  and $\tilde Y_0\deff \ti D_0\backslash \tilde \calr_0$ so that $(\tilde Y_0,{\c}_0, Y_0)$ is a finite regular cover.
Fix $z_*\in Y_0$ and let $w_*\in \tilde Y_0$ be a fixed preimage of $z_*$ by $\c_0$.
For every path $\alpha$ starting at $z_*$ we let $\c _0^* \alpha$ denote the lifted path of $\alpha$ starting at $w_*$.
Let $\isl_*:\pi_1(Y_0,z_*)\to \pi_1(Y_1,z_*) $ be the canonical morphism and $K\deff \{[\gamma]\in \pi_1(Y_0,z_*) \tq \c _0^*\gamma \mbox{ is closed in }\tilde Y_0
\}$ be the subgroup defined in \refrema{corresp_cover}.

\begin{defi}
\label{equiv_path}
Let $\alpha^1_{z_*z}$, $\alpha^2_{z_*z}$ be two paths from $z_*$ to $z\in Y_1$. 
We say that $\alpha^1_{z_*z}\sim\alpha^2_{z_*z}$ if $\alpha^2_{z_*z}\cdot\left(\alpha^{1}_{z_*z}\right)^{-1}\in \isl_*(K)$.

\end{defi}
In other words two paths $\alpha^1_{z_*z}$ and $\alpha^2_{z_*z}$ based at $z_*$ are equivalent if $\alpha^2_{z_*z}\cdot\left(\alpha^{1}_{z_*z}\right)^{-1}$ is homotopic inside $Y_1$ to a loop $\gamma \subset Y_0$ such that its lifted path based at some preimage of $z_*$ is closed.
Note that in general a loop $\gamma \subset Y_0$ starting at $z_*$ may be homotopic inside $Y_1$ to the constant path while its lifted path by $\c_0$ based at any preimage of $z_*$ is not closed.
This definition gives an equivalence relation between paths in $Y_1$.
{The equivalent class of a path $\alpha_{z_*z}$ in $Y_1$ starting at $z_*$ is denoted by $[\alpha_{z_*z}]$.}
If $\alpha^1_{z_*z}$ and $\alpha^2_{z_*z}$ are homotopic inside $Y_1$ then they are equivalent in the sense of \refdefi{equiv_path}.
Set $\tilde Y_1 \deff \{ [\alpha_{z_*z}] \tq z\in Y_1  \}$ and $\c_1 ([\alpha_{z_*z}])= z$.
The map $\c_1 : \tilde Y_1 \to Y_1$ is obviously well-defined.

We define a topology on $\tilde Y_1$ as follows. 
A subset $\tilde U$ of $\tilde Y_1$ is a neighbourhood of $[\alpha_{z_*z}]$ if there exists a contractible open set $\calu$ containing $z$ in $Y_1$ such that
\[
 \{[\alpha_{z_*z}\cdot\beta] \tq \beta \mbox{ is a path in } \calu \mbox{ starting at }z  \} \subset \tilde U.
\]
\noindent $\tilde Y_1$ becomes a connected topological space and the map $\c_1$ is open and continuous.

Let us prove that $\tilde Y_1$ is Hausdorff. Let $[\alpha^1_{z_*z_1}] \neq [\alpha^2_{z_*z_2}]$. 
We shall prove there exist two disjoints neighbourhoods $\tilde U_1$ of $[\alpha^1_{z_*z_1}]$ and $\tilde U_2$ of $[\alpha^2_{z_*z_2}]$ in $\tilde Y_1$. 
The case $z_1 \neq z_2$ is obvious so we can suppose that $z_1=z_2=z$.
Let $\calu$ be an open contractible neighbourhood of $z$ and suppose that there exist $\tilde z \in \calu$ and two paths $\beta^1_{z\tilde z},\beta^2_{z \tilde z}$ between $z$ and $\tilde z$ in $\calu$ such that
\begin{equation}\label{hausd}
 \left[ \alpha^1_{z_* z}\cdot \beta^1_{z\tilde z} \right] = \left[ \alpha^2_{z_* z} \cdot\beta^2_{z\tilde z}\right].
\end{equation}
The following equality holds in $\pi_1(Y_1,z_*)$
\[\left[ \alpha^2_{z_* z}\cdot \beta^2_{z\tilde z}\cdot \left(\alpha^1_{z_* z}\cdot \beta^1_{z\tilde z}\right)^{-1}  \right]_{\pi_1(Y_1,z_*)}=  \left[ \alpha^2_{z_* z}\cdot  \left(\alpha^1_{z_* z}\right)^{-1}  \right]_{\pi_1(Y_1,z_*)}\]
and \eqref{hausd} implies that $\alpha^2_{z_* z} \cdot \left(\alpha^1_{z_* z}\right)^{-1}\in \isl_*(K)$.
Hence $[\alpha^1_{z_*z_1}] = [\alpha^2_{z_*z_2}]$ which gives us a contradiction.
Then $\tilde U_1 \deff \{[\alpha^1_{z_*z}\cdot\beta] \tq \beta \subset \calu
\}$ and $\tilde U_2 \deff \{[\alpha^2_{z_*z}\cdot\beta] \tq \beta \subset \calu
\}$ are the desired open sets. $\tilde Y_1$ is a Hausdorff topological space.

Let us prove that $(\tilde Y_1,\c _1,Y_1)$ is a regular topological cover. 
Fix $z\in Y_1$ and let $\left\{[\alpha_{z_*z}^i]\right\}_{i\in I}$ be the preimages set of $z$ by $\c _1$ where $I$ is at most countable.
Let $\calu$ be a contractible open neighbourhood of $z$ in $Y_1$.
One obtains that 
\[\c _1^{-1}(\calu)=\bigcup_{i\in I} \{ [\alpha_{z_*z}^i\cdot \beta] \tq \beta \mbox{ is a path in } \calu \mbox{ starting at }z \}.\]
The fact that $\calu$ is contractible implies that the previous union is disjoint.
Therefore $\c _1$ is a regular cover and $\tilde Y_1$ inherits a structure of complex manifold.

\smallskip
For every $w\in \tilde Y_0$ we set $\tilde \isl(w)\deff[\c_0\circ \lambda_{w_*w}]$ where $\lambda_{w_*w}$ is a path joining $w_*$ and $w$ in $\tilde Y_0$.
Then by construction of the equivalence relation the map $\tilde \isl:\tilde Y_0 \to \tilde Y_1$ is well-defined, continuous and verifies $\c_1 \circ \tilde\isl=\c_0$.
Let us prove that its restriction on the fiber $\c _0^{-1}(z_*)$ is surjective onto $\c _1^{-1}(z_*)$.
Let $[\alpha]$ be a preimage of $z_*$ by $\c _1$.
Since $\isl_*:\pi_1(Y_0,z_*)\to \pi_1(Y_1,z_*)$ is surjective the loop $\alpha$ is homotopic inside $\pi_1(Y_1,z_*)$ to $\gamma \subset Y_0$.
It follows that $\alpha$ and $\gamma$ are equivalent in the sense of \refdefi{equiv_path}.
Denoting $w\deff \c _0^* \gamma (1)$ one deduces that $\tilde{\isl}(w)=[\gamma]=[\alpha]$ and $\tilde{\isl}: \c _0^{-1}(z_*) \to \c _1^{-1}(z_*)$ is surjective.
By using the lifting property of regular covers one obtains that $\tilde \isl:\tilde Y_0 \to \c_1^{-1}(Y_0)$ is surjective.
In this case $(\tilde Y_1,\c _1, Y_1)$ is a finite regular cover with $b_1$ sheets {where $b_1$ denotes the cardinal of $\c_1\inv(z)$ which doesn't exceed $b_0$.}
$\deff \# I \leq b_0$.
By \refthm{ext_norm_stein} and {\sl Remark \ref{grauert_connected}} it uniquely extends to $b_1$-sheeted connected analytic cover $(\tilde D_1,\c _1, D_1)$ over $D_1$. 
Let us prove that $\tilde \isl:\tilde Y_0 \to \tilde Y_1$ can be extended to $\tilde D_0$. 
Fix $w_0\in \tilde{\calr_0}$ and set $z_0\deff \c_0(w_0)$.
Let $\tilde U$ be an open connected neighbourhood of $w_0$ inside $\tilde D_0$ such that $\tilde U \backslash \tilde{\calr_0}$ is connected and let $\tilde V$ be a connected component of $\c _1^{-1}(U)$ inside $\tilde D_1$ such that
$\tilde{\isl}:\tilde U \backslash \tilde{\calr}_0 \to \tilde V$ is holomorphic, where $U\deff \c _0(\ti U)$.
After shrinking $\tilde U$ we can suppose that $\tilde V$ is biholomorphic to an analytic subset $A$ of some bounded domain $\Omega\subset \cc^m$.
Then one can replace $\tilde{\isl}$ by $m$ bounded holomorphic maps $\tilde{\isl}_k:\tilde U \backslash\tilde{\calr}_0 \to \cc$ ($k=1,\dots,m$).
By normality of $\tilde U \subset \tilde{D}_0$ every $\tilde{\isl}_k$ can be extended to a holomorphic function $\underline{\tilde{\isl}}_k:\tilde U\to \cc$.
The map $(\underline{\tilde{\isl}}_1,\dots,\underline{\tilde{\isl}}_m):\tilde U \to A$ extends $\tilde{\isl}:\tilde Y_0 \to \tilde Y_1$ on the whole of $\tilde D_0$.
Denoting by $\tilde \isl$ that extension one obtains that $\tilde{\isl}:\tilde D_0 \to \tilde D_1$ is holomorphic and diagram \ref{diag_D} is commutative. The first statement of \refthm{main_thm} is proved.

Let us prove statement \emph{\ref{b_1=b_0}}. 
If $b_1=b_0$ the restriction of the map $\tilde \isl:\tilde Y_0 \to \tilde Y_1$ on the fiber $\c _0^{-1}(z_*)$ becomes bijective onto $\c _1^{-1}(z_*)$.
It follows that $\tilde \isl:\tilde Y_0 \to \tilde Y_1$ is globally injective and $(\tilde D_1,\c _1, D_1)$ is an extension in the strong sense of $(\tilde D_0,\c _0, D_0)$ over $D_1$.

\smallskip

Let $(\tilde D'_1,\c' _1, D_1)$ be another {connected} extension of $(\tilde D_0,\c _0, Y_0)$. 
We shall proved that the number of its sheets can not be larger than $b_1$ \ie $(\tilde D_1,\c _1, D_1)$ is the unique maximal analytic cover which extends $(\tilde D_0,\c _0, D_0)$.
There exists a holomorphic map $\tilde \isl' : \tilde D_0 \to \tilde D'_1$ such that the following diagram
\begin{equation*}
 \xymatrix{
    \tilde D'_1 \ar[d]_{\c'_1} & \tilde D_0\ar[r]^{\tilde\isl} \ar[l]_{\tilde\isl'} \ar[d]^{\c_0} & \tilde D_1 \ar[d]^{\c_1} \\
    D_1 & D_0 \ar@{^{(}->}[r]_{\isl} \ar@{_{(}->}[l]^{\isl}& D_1
  }
\end{equation*}
is commutative \ie $\c'_1\circ\tilde\isl'=\c _0=\c_1\circ\tilde\isl$.
Set  $\tilde Y'_1\deff \c_1^{'-1}(Y_1)$. The restriction  $\c'_1=\restriction{\c' _1}{\tilde Y'_1}:\tilde Y'_1 \to Y_1$ is a regular analytic cover over $Y_1$.
Let $\zeta'_*$ be a preimage of $z_*$ by $\c _1'$.
We shall prove the following 
\begin{lem}\label{closed}
Let $\alpha_1$, $\alpha_2$ be two equivalent paths in the sense of \refdefi{equiv_path}.
Then its lifted path $\c _1'^{*}\left(\alpha_2\cdot \alpha_1^{-1} \right)$ starting at $\zeta'_*$ is closed.
\end{lem}

{\slsf Proof of \lemma{closed}.}
Let $\gamma \subset Y_0$ be a loop such that $\alpha_2\cdot\alpha_1^{-1}$ is homotopic inside $Y_1$ to $\gamma$ and such that $\c_0^*\gamma$ is closed.
Since $\c '_1\circ\tilde{\isl}'=\c _0$ it follows that $\c^{'*}_1\gamma$ is closed
and one deduces that $\c _1'^{*}\left(\alpha_2\cdot \alpha_1^{-1} \right)$ is closed.
\lemma{closed} is proved.

\qed

\smallskip
Now let us prove statements \emph{\ref{b'_1<b_1}} and \emph{\ref{b'_1=b_1}}.
Let $\alpha \subset Y_1$ be a path starting at $z_*$. 
By the previous Lemma the map
\[ \fonction{\Psi}{\tilde Y_1}{\tilde Y'_1}{[\alpha]}{\c _1'^{*}\alpha (1)}
\]
is well-defined, continuous, surjective and it verifies $\c'_1\circ \Psi=\c _1 $.
Hence $(\tilde Y'_1,\c '_1,Y_1)$ is a finite regular cover and the number $b'_1$ of its sheets can not be larger than $b_1$.
Remark that if $b'_1=b_1$ the restriction $\Psi:\c_1 ^{-1}(z)\to \c_1 ^{'-1}(z)$ on every fiber is a surjective map between finite sets with $b_1$ elements. It is injective and $\Psi: {\tilde Y_1}\to{\tilde Y'_1}$ becomes bijective.
Its inverse is holomorphic given by $\Psi^{-1}:\zeta\in\tilde{Y}'_1 \mapsto \left[\c' _1\circ\tilde{\alpha}_{\zeta'_*\zeta}\right]$ where $\tilde{\alpha}_{\zeta'_*\zeta}$ is path between $\zeta'_*$ and $\zeta$ inside $\tilde{Y}'_1$.
One deduces that $\Psi$ is a biholomorphism which uniquely extends to a biholomorphic map $\Psi:\tilde D_1 \to \tilde D'_1$ such that $\c '_1\circ\Phi=\c_1$.

\qed

\smallskip

\begin{rema} \rm
In the assumptions of \refthm{main_thm} if the analytic cover $(\ti D_0,\c_0,D_0)$ is Galois then $(\ti D_1,\c_1,D_1)$ is Galois.
\end{rema}

The following example is well known.
It shows that an analytic cover can be extended only by gluing the sheets.
Remark that the branching divisor in it is empty.

\begin{exmp} \rm
\label{gluying} 
Take as $X=\cc^2\setminus \rr^2$ and as $\tilde X$ a $b$-sheeted regular analytic cover of $X$. It cannot be extended over any point of $\rr^2$ in the strong sense because holomorphic functions on
this $\tilde X$ do not separate points. But it obviously extends to
a trivial cover $(\cc^2,\id , \cc^2)$ after ``gluing the sheets".
\end{exmp}

\smallskip
%

\newsect[POLYDISK]{Extension of analytic covers over a polydisk}

\newprg[POLYDISK_PIC_SIM]{Lemma of Picard-Simart}

The following statement is due to \cite{PS}. It is proved in \cite{Ni}.

\begin{lem} \rm \label{pic_sim_0}
Let $\Delta^n \subset \cc^n$ be the unit polydisk and $\Gamma \deff \{|w|<1 \} \subset \cc$ be the unit disk.
Let $\calr$ be a one-codimensional analytic subset of $\Delta^n\times \Gamma$ such that $\calr \cap ( \Delta^{n}\times \d\Gamma)=\emptyset$.
Then there exists a proper analytic subset $\sigma \subset \Delta^n$ such that for every $z_*\in \Delta^n \backslash \sigma$ any one-dimensional closed curve $\gamma$ in $\left(\Delta^n\times \Gamma\right) \backslash\calr$ can be continuously deformed to a closed curve in $\left(\{z_*\}\times \Gamma \right)\backslash\calr$.
\end{lem}

By using the methods of the proof of \lemma{pic_sim_0} given by \cite{Ni} one deduces a more precise statement about surjectivity of the natural map between fundamental groups.

\begin{thm}
\label{pic_sim}
In the assumptions of \lemma{pic_sim_0}
 there exists a nowhere dense subset $A \subset \Delta^n$ such that for every $z_*\in \Delta^{n}\backslash A$ and every $Z_*=(z_*,w_* )\in \left(\Delta^n \times \Gamma\right)\backslash \calr$ the natural homomorphism
\[\isl_*:\pi_1\left[(\{z_* \}\times\Gamma) \backslash\calr,Z_*\right]\to \pi_1\left[\left(\Delta^n \times \Gamma\right)\backslash\calr,Z_*\right]\]
is surjective \ie every loop $\gamma$ in $\left(\Delta^n \times \Gamma\right)\backslash \calr$ starting at $Z_*$ is homotopic to some path $\gamma^* \subset (\{z_* \}\times\Gamma) \backslash\calr$ within the loops starting at $Z_*$.
\end{thm}

\proof 
We denote the standard coordinates in $\cc^{n+1}$ as $z=(z_1,\dots,z_n)=(x_1+iy_1,\dots,x_n+iy_n)\in \Delta^n$ and $w=u+iv \in \Gamma$.
By the hypothesis on $\calr$ there exists a monic Weierstrass polynomial 
\[P(z,w)=w^{\nu}+\sum_{i=1}^{\nu}a_i(z)w^{\nu-i}\]
 such that $a_i$ is holomorphic on $\Delta^n$, $P(z,w)$ has no multiple factor and
  \[\calr=\{(z,w)\in \Delta^n\times\Gamma \tq P(z,w)=0
\}.\]
Since $\Delta$ is biholomorphic to the square $\{z=x+iy \in \cc \tq -1<x<1 \text{ and } -1<y<1  \}$ we may assume that $\Delta^n$ is the cube.
Set $\sigma \deff \left\{ z\in \Delta^n \tq \discr_w P(z,w)=0  \right\}$. 
The following statement  is proved in \cite{Ni} as Lemma 2.9.
\begin{lem} \rm \label{Ni_2.9}
Let $D$ be a domain of $\cc^{n+1}$ whose coordinates are denoted as $z=(z_1,\dots,z_n)\in\cc^n$ and $w\in \cc$.
Let $\calr=\{f(z,w)=0\}$ be an analytic hypersurface of $D$.
There exists a linear transformation $\phi$ of $\cc^{n+1}$ such that in the coordinates $(z',w')=\phi(z,w)$ one has
\[ \forall (a'_1,\dots,a'_n,b')\in \cc^{n+1} \; \left[ f\circ\phi^{-1}(a'_1,\dots,a'_n,b')=0 \Rightarrow f\circ\phi^{-1}(a'_1,\dots,a'_n,w')\not\equiv 0  \right].
\]
\end{lem}

We prove \refthm{pic_sim} by induction.

\smallskip

\noindent{\slsf Case $n=1$.}
According to \lemma{Ni_2.9} we can suppose after taking a linear transformation that $\calr$ does not contain any complex hyperplane of the form $w=d$ where $d\in \Gamma$ is constant.
The set $\sigma \subset \Delta$ consists of a countable number of points $A_k=A'_k+iA''_k, k\in \nn$.
Set $A\deff\{x+\isl y \in \Delta \tq \exists k\in\nn  \;  x=A'_k\text{ or }y=A''_k  \}$.
Let $Z_*=(x_*, y_*,u_*, v_*)\in \left(\Delta \times \Gamma\right)\backslash \calr$ be such that $x_*+\isl y_*\notin A$ and $\gamma$ be a loop in $\left(\Delta^n \times \Gamma\right)\backslash \calr$ starting at $Z_*$.
We can suppose that $\gamma$ is a real analytic path and its projection onto every real axis $x$, $y$, $u$ and $v$ is not reduced to a single point.
For every $M=(x_M,y_M,u_M,v_M)\in \gamma$ we set
$X_M\deff\{y=y_M\}\cap \Delta$ and we let $\calx_M$ denote the cylinder $ X_M\times \Gamma$.
We have the following.
\begin{lem}\label{lem_segment}
For every $M=(x_M,y_M,u_M,v_M)\in \gamma$ we can find a real one-dimensional open line segment $L(M)\subset\calx_M$ containing $M$ 
of the form
\begin{equation}\label{segment}
L(M)\deff \left\{\left(x,y_M,u_M+\alpha_M(x-x_M),v_M+\beta_M(x-x_M)\right)\tq x\in[-1,1]
\right\}
\end{equation}
where $\alpha_M$ and $\beta_M$ are real constants and such that :
\begin{enumerate}[label=\emph{(\arabic*)}]
\item $L(M)\cap \calr = \emptyset$ and \label{i}
\item $L(M)\cap \left[\{|x|\leq 1
\}\times \{y_M\}\times \d\Gamma
\right]=\emptyset.$ \label{ii}
\end{enumerate}
\end{lem}

Let $M_0\in \gamma$ and let $\alpha_0$ and $\beta_0$ be the real constants in the definition of $L(M_0)$.
There exists a subarc $[M_0'M_0'']$ of $\gamma$ which contains $M_0$ as an interior point and such that for every $M\in[M_0'M_0'']$ 
we can take $\alpha_M=\alpha_0$ and $\beta_M=\beta_0$ in the definition \eqref{segment}   of $L(M)$.
By compactness of $\gamma$ we can find a finite number of points $M_0,M_1,\dots,M_q$ such that:
\begin{itemize}
\item $\gamma=\bigcup_{i=0}^{q-1} [M_iM_{i+1}]$,

\item $M_0=M_{q}=Z_*$ and

\item for every point $M \in [M_i M_{i+1}]$ we can take $\alpha_M=\alpha_i$ and $\beta_M=\beta_i$ where $\alpha_i,\beta_i$ correspond to the definition of $L(M_i)$ in \eqref{segment}.
\end{itemize}

\begin{figure}
\centering
\includegraphics[width=0.8\linewidth]{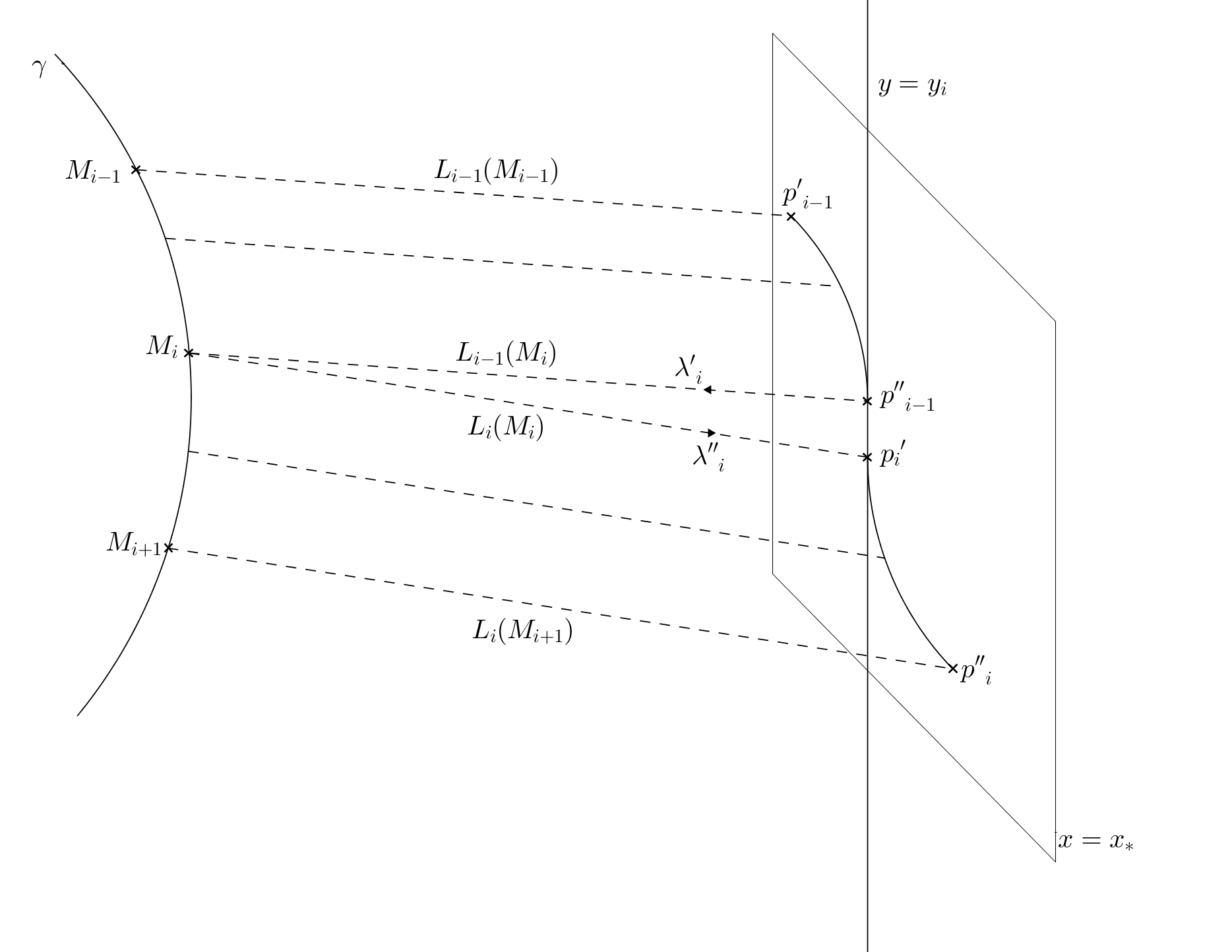}
\caption{}
\label{proj}
\end{figure}

If $M\in[M_i,M_{i+1}]$ the line segment $L(M)$ is denoted by $L_i(M)$.
We denote $M_i=(x_i, y_i,u_i, v_i)$.
Since the projection of $\gamma$ on the axes $x$ and $y$ does not reduce to a point we can suppose after a small perturbation that $M_i \notin A\times\Gamma$ for $i=1,\dots,q$.
For every $M=(x_M, y_M,u_M , v_M)\in [M_iM_{i+1}]$ we define its projection to the real hyperplane $\{x=x_*\}$
\[p_i(M)\deff \restriction{L_i(M)}{x=x_*}=(x_*,y_M,u_M+\alpha_i(x_*-x_M),v_M+\beta_i(x_*-x_M)).\]
Set $p_i'\deff p_i(M_i)$ and $p_i''\deff p_i(M_{i+1})$.
Remark that every $M_i$ has two projections $p'_i$ and $p''_{i-1}$ in $\{x_*,y_i\}\times\Gamma$. 
We can continuously deform the arc $[M_iM_{i+1}]$ to the arc $[p_i'p_i'']\deff \{p_i(M) \tq M \in [M_i,M_{i+1}]  \}$ in such a way that $M\in[M_iM_{i+1}]$ moves to $p_i(M)$ along $L_i(M)$, see Figure \ref{proj}.
By hypothesis on $L_i(M)$ this deformation is in $\left(\Delta \times \Gamma\right)\backslash \calr$.
Since $p_0'=p_{q-1}''=Z_*$ the arc $[M_0 M_1]$ can be deformed to $[M_0p_0'']$ within the paths starting at $M_0$.
Likewise $[M_{q-1} M_0]$ can be deformed to $[p'_{q-1}M_0]$ within the paths ending at $M_0$.
Let $\lambda'_i$ be the path between $p_{i-1}''$ and $M_i$ on the segment $L_{i-1}(M_i)$ and $\lambda''_i$ a path between $M_i$ and $p_i'$ on $L_i(M_i) $.
Then denoting $\lambda_i\deff \lambda_i'\cdot\lambda''_i\subset \calx_{M_i}\backslash\calr$ one obtains that $\gamma$ is homotopic in $\left(\Delta \times \Gamma\right)\backslash \calr$ to the loop
\[ \hat{\gamma}\deff[M_0p_0'']\cdot \lambda_1 \cdot[p'_1p_1'']\cdots \lambda_{q-1}\cdot [p'_{q-1} M_0]
\]
within the loops starting at $M_0=Z_*$.

We have the following.

\begin{lem}\label{assert_I}

Let $D\deff ]-1,1[\times \Gamma \subset \rr^3$ be a cylinder and $\call_j: (t,u_j(t),v_j(t))$ a smooth arc in $\overline{D}$ ($j=1, \dots, \nu$).
We make the following assumptions.
\begin{enumerate}[label=\emph{(\arabic*)}]
\item $\call_j \cap ([-1,1]\times\d\Gamma)=\emptyset$ for every $j=1,\dots,\nu$ and
\item $\call_j\cap\call_k=\emptyset$ if $j\neq k$.
\end{enumerate}
Set $\call\deff\bigcup_{j=1}^{\nu}\call_j$ and let  $\gamma$ be a path in $D\backslash\call$ such that $\gamma(0),\gamma(1)\in \{t_0
\}\times\Gamma$ for some $t_0 \in ]-1,1[$.
Then $\gamma$ is homotopic to a path $\tilde \gamma \subset \{t_0\} \times \Gamma$ within the paths starting at $\gamma(0)$ and ending at $\gamma(1)$.
\end{lem}

The proof of this Lemma consists in observing that there exists a homeomorphism $\Phi:D\to D$ such that
\begin{itemize}
\item $\forall t\in ]-1,1[ \; \Phi(\{t\}\times \Gamma)=\{t\}\times \Gamma$,
\item $\Phi(\call_j)=]-1,1[\times \{(a_j,b_j) \}$ where $(a_j,b_j)\neq (a_k,b_k)$ if $j\neq k$ and
\item the restriction $\Phi:{\{t_0\}\times\Gamma }\to {\{t_0\}\times\Gamma }$ is the identity mapping.
\end{itemize}

On deduce that we can deform $\hat{\gamma}$ in $\left(\Delta \times \Gamma\right)\backslash \calr$ to a closed curve $\tilde{\gamma}$ in $ \caly(z_*)\backslash\calr$
where 
\begin{equation*}
\caly(z_*)\deff \left(\Delta\times\Gamma\right) \cap\{x=x_*\}.
\end{equation*}
Indeed since $(x_i, y_i) \notin A$ we have $\calx_{M_i}\cap \sigma = \emptyset$.
Then $\calr\cap \calx_{M_i}$ consists of an union of $\nu$ real analytic arcs $\call_j$ verifying the assumptions of \lemma{assert_I}.
Moreover for every $i\in\{1,\dots,q\}$ $\lambda_i(0)=p''_{i-1}\in \{x_*,y_i\}\times\Gamma $ and $\lambda_i(1)=p'_{i}\in \{x_*,y_i\}\times\Gamma $.
Since $\lambda_i$ does not intersect $\calr$ \lemma{assert_I} implies that for every $i=1,\dots,q-1$ there exists a path $\tilde{\lambda_i}\subset (\{x_*, y_i\}\times\Gamma)\backslash\calr $ starting at $p_{i-1}''$ and ending at $p_i'$ such that $\lambda_i$ is homotopic to $\tilde\lambda_i$ inside $\calx_{M_i}\backslash\calr$ with the same initial and terminal points.
Thus $\gamma$ is homotopic in $\left(\Delta \times \Gamma\right)\backslash \calr$ to
\[
\tilde{\gamma}\deff[M_0p_0'']\cdot \tilde\lambda_1 \cdot[p'_1p_1'']\cdots \tilde\lambda_{q-1}\cdot [p'_{q-1} M_0] \subset \caly(z_*)\backslash\calr.
\]
The last thing is to prove that $\tilde{\gamma}$ can be deformed to a close path $\gamma^*$ in $\{x_*,y_*\}\times\Gamma$.
But this is an application of \lemma{assert_I} to $\tilde \gamma$ by taking $t=y$.
\refthm{pic_sim} is proved when $n=1$.

\smallskip

\noindent{\slsf Case $n\geq 2$}.
Let $z_*=(z_1^*,\dots,z_n^*) \in \Delta^n\backslash\sigma$ where $z_n^*=x_n^* +iy_n^*$.
According to \lemma{Ni_2.9} we can suppose after taking linear transformation that $\calr$ does not contain any complex lines of the form $\{z_1=c_1,\dots,z_{n-1}=c_{n-1}, w=d\}$ where $(c_1,\dots,c_{n-1})\in\Delta^{n-1}$ and $d\in \Gamma$ are constants.
Let $\gamma$ be a closed curve starting at $Z_*=(z_*,w_*)\in \left(\Delta^n \times \Gamma\right)\backslash \calr$.
As previously we can suppose that $\gamma$ is real analytic of the form $\gamma(t)= (\varphi(t),\chi'(t),\chi''(t),\psi'(t),\psi''(t)) $
{where $\varphi:  [0,1] \to \Delta^{n-1}$ and $\chi',\chi'',\psi',\psi'':[0,1] \to ]-1,1[$} are real analytic functions such that the projection of $\gamma$ onto every axis $x_1,y_1,\dots,x_n,y_n,u$ and $v$ is not reduced to a single point.
We shall prove there exists a nowhere dense subset $A$ of $\Delta^n$ such that the assumption of the lemma is satisfied if $z_*\notin A$.

Fix $M \in \gamma$ and denote $M=(z'_M,x_M,y_M,u_M,v_M)=(\varphi(s),\chi'(s),\chi''(s),\psi'(s),\psi''(s))$.
Set $X_M\deff \{(\varphi(s),x,\chi''(s))\in \Delta^n \tq -1\leq x \leq 1 \}$ and $\calx_M\deff X_M \times \Gamma$.
There exists a real one-dimensional segment $L(M) \subset \calx_M $ containing $M$ of the form
\begin{equation}
L(M)\deff \{(z'_M,x,y_M,u_M+\alpha_M(x-x_M),v_M+\beta_M(x-x_M)\tq -1\leq x \leq 1 \}
\end{equation} 
where $\alpha_M$ and $\beta_M$ are real constants such that 
\begin{enumerate}
\item $L(M)\cap \calr=\emptyset$ and
\item $L(M)\cap \left[  \{z_M' \} \times \{ |x|\leq 1 \} \times \{y_M\}\times \d\Gamma \right]=\emptyset $.  
\end{enumerate}

Set $\Lambda_{x_n^*}\deff \Delta^{n-1}\times \{x_n^*\}\times \{|y_n|\leq 1 \}\times \Gamma$, $\caly_M(z_*)\deff \{ (\varphi(s),x_n^* , y)\in \Delta^n \tq -1 \leq y \leq 1  \}\times \Gamma$ and $\tilde{\caly}_{\gamma}=\bigcup_{M\in\gamma}\caly_M(z_*)$.
According to the method used in the proof of the case $n=1$ there exists a nowhere dense subset $A'\subset \{ |x_n|<1, |y_n|< 1  \}$ such that for every $(x_n^*,y_n^* )\notin A'$ the curve $\gamma$ is homotopic to a closed path $\tilde{\gamma}:t\mapsto (\varphi(t),x_n^*,y(t),u(t),v(t))$ in $\tilde{\caly}_{\gamma}\backslash\calr$ where $y(t),u(t)$ and $v(t)$ are continuous.

Then $\tilde{\gamma}$ can be deformed in $\Lambda_{x_n^*}\backslash\calr$ to a curve $\hat{\gamma}$ in $\left( \bigcup_{M\in\gamma} \{ (\varphi(t),x_n^*,y_n^*) \}\times\Gamma \right)\backslash\calr$.
Denoting $\Lambda^{n-1}\deff \Delta^{n-1}\times\{z_n^*\}\times \Gamma$ and $\calr^{n-1}\deff\calr\cap\Lambda^{n-1}$ one obtains that $\hat{\gamma}\subset \Lambda^{n-1}\backslash\calr^{n-1}$ and the proof of the Lemma reduces to the case $n-1$.
By induction \refthm{pic_sim} is proved.

\qed

\newprg[POLYDISK_result]{Proof of Theorem \ref{ext_disk}}

We prove in this section \refthm{ext_disk} from Introduction.

According to \refthm{pic_sim} there exists $z_*\in D_0\backslash\calr$ such that the natural morphism $\isl_*:\pi_1(D_0\backslash\calr,z_*)\to \pi_1(\Delta^n\backslash\calr,z_*) $ is surjective.
By \refthm{main_thm} one obtains \refthm{ext_disk_b}.
The second part of the Theorem is a particular case of the following statement, see \refthm{thm_ext_2}

\begin{thm}\label{thm_cover_2}
Let $(\tilde D,\c , D)$ and $(\tilde D',\c',D)$ be $2$-sheeted {connected} analytic covers over a domain $D\subset\cc^n$ with the same ramification divisor $\calr \subset D$.
Then the covers are equivalent \ie there exists a biholomorphic map $\tilde{\isl}:\ti D \to \ti D'$ such that $\c '\circ\tilde{\isl}=\c$.
\end{thm}

\proof The proof consists from the following steps.

\smallskip

\noindent{\slsf Step 1 :} Let us prove the result when $D=\Delta^n$ is the unit polydisk and $\calr$ is such that $\calr\cap\left(\Delta^{n-1}\times\d\Delta\right) = \emptyset$.
The divisor $\calr$ is the zero set of a Weierstrass polynomial of degree $\nu$ with respect to $z_n$.
According to \refthm{pic_sim} there exists $z'_*\in \Delta^{n-1}$ such that the natural morphism
$\isl_*:\pi_1(\Delta_{z_*'}\backslash\calr,z_*)\to \pi_1(D\backslash\calr,z_*)$ is surjective {where $z_*=(z_*',z_n)\in D\backslash\calr$ and $\Delta_{z_*'}\deff {z_*'}\times \Delta$.}
Set $\tilde{Y}\deff \c ^{-1}(D\backslash\calr) $ and $\tilde{Y}'\deff \c ^{'-1}(D\backslash\calr)$.
Let $w_*\in \tilde Y$ and $w'_*$ be some preimages of $z_*$ by $\c$ and $\c'$ respectively. 
The point $z'_*$ has been chosen so that $\Delta_ {z'_*}\cap \calr$ is a finite set of $\nu$ elements $\left\{a_1,\dots,a_{\nu} \right\}$.
Fix $i\in \{1,\dots,\nu\}$ and take a neighbourhood $U_i$ of $a_i$ inside $\Delta_{z'_*}$ such that $\c^{-1}(U_i\backslash \calr)$ is connected.
Let $z_i\in U_i\setminus\calr$.
The generator $[\alpha_i]$ of $\pi_1(U_i\backslash\calr,z_i)\simeq \zz$ is such that
the lift of $\alpha_i$ at any preimage of $z_i$ by $\c$ and $\c'$ is not closed.
Let $\lambda_i$ be a path between $z_*$ and $z_i$ inside $\Delta_{z'_*}\backslash\calr$ and define $\gamma_i \deff \lambda_i\cdot\alpha_i\cdot \lambda_i^{-1}$.
By construction the lift $\tilde{\gamma}_i$ of $\gamma_i$ starting at $w_*$ by $\c$ is not closed.
By using the same argument one obtains that the lift $\tilde{\gamma}'_i$ by $\c'$ starting at $w'_*$ is not closed.
Moreover $\pi_1(\Delta_{z_*'}\backslash\calr,z_*)$ is the free group generated by $\{[\gamma_1],\dots,[\gamma_{\nu}]\}$.
It is easy to prove by induction that the lift of any loop $\gamma_{i_1}^{n_1}\cdots \gamma_{i_s}^{n_s} \subset \Delta_{z_*'}\backslash\calr$ by $\c$ is closed if and only if $n_1+\dots +n_s$ is even.

One deduces that if two paths $\gamma_{i_1}^{n_1}\cdots \gamma_{i_s}^{n_s}$ and $\gamma_{j_1}^{m_1}\cdots \gamma_{j_q}^{m_{s'}}$ are homotopic inside $D\backslash\calr$ then $n_1 +\dots +n_s -(m_1 + \dots + m_{s'})$ is even.
Moreover \refthm{pic_sim} implies that every loop $\gamma$ in $D\backslash\calr$ starting at $z_*$  is homotopic inside $D\backslash\calr$ to $\gamma_{i_1}^{n_1}\cdots \gamma_{i_s}^{n_s}\subset \Delta_{z_*'}\backslash\calr$ for some $n_j \in \zz^*$ ($j=1,\dots,s$).
According to the previous remark  one deduces that the parity of $n_1 + \dots + n_s$ does not depend of the choice of such a decomposition.
Then
$\c^*\alpha$ is closed if and only if $n_1+\dots + n_s$ is even.
In particular the lift $\c^*\alpha$ of any loop $\alpha \subset D\backslash \calr$ by $\c$ is closed if and only if $\c^{'*}\alpha$ is closed.
The regular covers $(\tilde Y,\c , D\backslash\calr)$ and $(\tilde Y',\c' , D\backslash\calr)$ are equivalent and according to \refthm{ext_norm_stein} and {\sl Remark \ref{grauert_connected}} the result follows.

\smallskip

\noindent{\slsf Step 2 : } We need to prove the following.

\begin{lem} 
\label{lem_exhaus_2}
{In the assumptions of \refthm{thm_cover_2} for every $p\in \calr$ there exist a neighborhood $U_p$ of $p$ inside $D$ and $z_p\in U_p\backslash\calr$ such that the fundamental group $\pi_1(U_p\setminus\calr,z_p)$ is generated by $\nu_p$ loops $\gamma_1,\dots,\gamma_{\nu_p}$ based at $z_p$.
Moreover for every $j=1,\dots,\nu_p$ and for any two preimages $w_p$ and $w'_p$ of $z_p$ by $\c $ and $\c'$ respectively the lifts $\c^*\gamma_j$, $\c^{'*}\gamma_j$ based at $w_p$ and $w'_p$ respectively are not closed.}

\end{lem}

\proof Fix $p\in \calr$.
One can find a complex line $L$ containing $p$ and a neighbourhood $\Delta$ of $p$ inside $L$ such that $\Delta \cap \calr = \{p\}$.
Take a complex orthogonal direction $L^\perp$ of $L$ at $p$ in $\cc^n$ and let $z=(z_1,\dots,z_n)$ be the coordinates chart corresponding to $L^\perp \times L$.
There exists $\eps >0$ such that $\Delta_{\eps}^{n-1}\times \Delta \subset D$ and $(\Delta^{n-1}_{\eps}\times \d \Delta)\cap\calr=\emptyset$.
By using the construction used in the previous step one deduces the result.
\smallskip \qed
Remark that the open sets $U_p$ in the lemma can be chosen so that $\{ U_p \tq p\in\calr\}$ is an open covering of $D$.

\smallskip

\noindent{\slsf Step 3 : Proof of \refthm{thm_cover_2}}
Fix $z_*\in D\backslash\calr$ and let $w_*,w'_*$ be two preimages of $z_*$ by $\c $ and $\c'$ respectively.
Let $\alpha$ be a loop inside $D\backslash\calr$ at $z_*$.
According to the previous step there exists a finite number of points $p_1,\dots,p_N$ in $\calr$ such that
$\alpha \subset \left(U_1\cup\dots\cup U_N\right)\backslash\calr$ where $U_k\deff U_{p_k}$ ($k=1,\dots,N$) is the open set defined in \lemma{lem_exhaus_2}.
One can find a decomposition $\alpha = \alpha_1\cdots\alpha_q$ where $\alpha_i$ is a path inside some $U_{k(i)}\setminus\calr$ where $1\leq k(i) \leq N$ such that $\alpha_i(1)=\alpha_{i+1}(0)$, $i=1,\dots,q-1$.
Fix $i\in \{1,\dots,q\}$ and take a point $z_i \in U_i\setminus\calr$ which verifies the assumptions of \lemma{lem_exhaus_2}.
Let $\lambda_i$ (respectively $\mu_i$) be a path inside $U_i\setminus\calr$ between $\alpha_i(0)$ and $z_i$ (respectively between $\alpha_i(1)$ and $z_i$).
According to the previous lemma the loop $\lambda_i\inv \cdot \alpha_i \cdot \mu_i$ based at $z_i$ is homotopic inside $U_i\setminus \calr$ to some $\gamma_{i,j_1}^{n_{i,1}}\cdots \gamma_{i,j_{s_i}}^{n_{i,s_i}}$ where $\gamma_{i,j_k}$ verifies the assumptions of the previous lemma and $n_{i,k}\in \zz$ ($1\leq k \leq s_i$).
Such a decomposition of $\lambda_i\inv \cdot \alpha_i \cdot \mu_i$ is not unique but one can prove like in Step 1 that the parity of $n_{i,1}+\cdots+n_{i,s_i}$ is constant.
Moreover the lift of $\lambda_i\inv \cdot \alpha_i \cdot \mu_i$ is closed if and only if $n_{i,1}+\cdots+n_{i,s_i}$ is even.
Finally one deduces that $\alpha$ is homotopic inside $D\backslash\calr$ to 
\begin{equation*}
\displaystyle{\lambda_1\left(\gamma_{1,j_1}^{n_{1,1}}\cdots \gamma_{1,j_{s_1}}^{n_{1,s_1}} \right)\mu_1\inv\lambda_2\left(\gamma_{2,j_1}^{n_{2,1}}\cdots \gamma_{2,j_{s_2}}^{n_{2,s_2}} \right)\cdots \mu_{q-1}\inv \lambda_q \left(\gamma_{q,j_1}^{n_{q,1}}\cdots \gamma_{q,j_{s_q}}^{n_{q,s_q}} \right)\mu_q\inv.}
\end{equation*}

By induction one can see that if that decomposition is homotopic to the constant path inside $D\setminus\calr$ then \begin{equation*}
S\deff \sum_{i=1}^q \sum_{k=1}^{s_i} n_{i,k}
\end{equation*}
is even.
Then the parity of this sum is independant of the choice of such a decomposition for $\alpha$.
Finally one can prove that the lift of $\alpha$ is closed if and only if $S$ is even.
\refthm{thm_cover_2} is proved.

\qed

\smallskip

\begin{thm}
\label{thm_ext_2}
Let $D_0$ be a domain and $(\tilde D_0,\c_0,D_0)$ be a $2$-sheeted connected analytic cover with ramification divisor $\calr_0\subset D_0$.
Suppose that $\calr_0$ extends to a globally defined hypersurface $\calr_1$ in $D_1\supset D_0$ \ie there exists a holomorphic function $f\in \calo(D_1)$ be such that $\calr_1 \deff \{ z\in D_1 \tq f(z)=0 \}$.
Then the cover $(\tilde D_0,\c_0,D_0)$ uniquely extends in the strong sense to a $2$-sheeted connected analytic cover over $D_1$ with ramification divisor $\calr_1$.

\end{thm}

\proof 
One can choose $f\in \calo(D_1)$ of multiplicity one such that $\tilde{\calr'}\deff \{(0,z)\in \cc \times D_1 \tq f(z)=0 \}$ does not locally separate $\tilde D_1'\deff \{ (\zeta,z)\in \cc\times D_1 \tq \zeta^2 = f(z) \}$.
Then $(\tilde D_1',\c ', D_1)$ is a $2$-sheeted analytic cover where $\c':\ti D_1' \to D_1$ is induced by the projection $\cc \times D_1 \to D_1$.
By \refthm{thm_cover_2} the restriction of this cover over $D_0$ is equivalent with $(\ti D_0,\c_0 ,D_0)$.
One deduces the result.
\qed

\smallskip
\refthm{ext_disk_2} is a direct consequence of the previous Theorem.

\newsect[HART]{Hartogs type extension of analytic covers}
Recall that a smooth real valued function $\rho$ in an open set $\Omega\subset \cc^n$ is called {\slsf $(n-q)$-convex} at point $z_* \in \Omega$ if its Levi form $L_{\rho}(z_*)$ has at least $q+1$ positive eigenvalues.
Let us prove \refthm{ext_hart} from Introduction.
It consists from the following steps.

\noindent{\slsf Step 1 : } 
Let us prove the following result.
\begin{lem}
\label{ext_q_concave}
Let $q\geq 2$ and $M=\{\rho = 0 \}$ be a smooth strongly $(n-q+1)$-convex hypersurface in a domain $D\subset \cc^n$.
Set $D^+=\{\rho>0\}$ and let $(\tilde D^+, \c_0, D^+)$ be a $b_0$-sheeted connected analytic cover over $D^+$. 
Suppose that the ramification divisor $\calr$ extends in a neighbourhood of every $p\in M$.
Then there exists a neighbourhood $U_p$ of $p$ such that the cover extends to a $b_1$-sheeted connected analytic cover over $D^+\cup U_p$, where $b_1\leq b_0$.
It verifies properties \ref{b_1<b_0 },\ref{b_1=b_0},\ref{b'_1<b_1} and \ref{b'_1=b_1} in \refthm{main_thm}.
\end{lem}

If $q \geq 3$ such a hypersurface $\calr$ always extends in $\Delta^n$, see Theorem 8.3 in \cite{ST}.

\proof
Let $V$ be the neighbourhood where $\calr$ extends.
We may assume that $p=0$ and every branch of $\calr$ in $V$ contains $0$.
Let $\Sigma$ be the complex tangent of $M$ at $0$.
There exists a non empty subspace $S\subset \Sigma$ on which the complex Hessian of $L_{\rho,0}$ is positive definite.
One can find a one-dimensional disk $\Delta\subset D^+$ centered at $0$ sufficiently close to $S$ such that $\Delta \cap \calr=\{0\}$. 
Take a complex orthogonal direction $L$ of $\Delta$ at $0$ in $\cc^n$ and let $z=(z_1,\dots,z_n)$ be the coordinates chart corresponding to $L\times\Delta$.
There exist $\eps >0$ such that $\Delta^{n-1}_{\eps } \times \Delta\subset D$ and $\calr \cap \left(\Delta^{n-1}_{\eps } \times\d \Delta\right)=\emptyset$.
Moreover one can find a polydisk $W\subset \Delta^{n-1}_{\eps }$ such that $W\times\Delta \subset D^+$.
By \refthm{ext_disk} the restriction of the cover $(\tilde D^+,\c_0 ,D^+) $ over $W\times\Delta$ can be extended to a $b_1$-sheeted connected analytic cover over $D^+ \cup \left(\Delta^ {n-1}_{\eps}\times\Delta\right)$ where $b_1\leq b_0$.

\smallskip
\qed

\smallskip 

\noindent{\slsf Step 2 : } One obtains the following result.
\begin{lem}\label{ext_vois}
In the statements of \lemma{ext_q_concave} $(\tilde D^+,\c _0,D^+)$ extends to a $b_{U}$-sheeted connected analytic cover of $D^+ \cup U$, where $U=\bigcup_{p\in M}U_p$ and $b_U\leq b_0$.
\end{lem}

\proof
Let $U_p$, $U_q$ be polydisks inside $D$ centered at $p \in M$ and $q \in M$ respectively such that  
\begin{enumerate}
\item $U_p\cap U_q \cap M \neq \emptyset$,

\item $(\tilde D^+, \c_0,D^+)$ extends to a $b_p$-sheeted cover $(\tilde U_p,\c _p, D^+ \cup U_p)$ and

\item $(\tilde D^+, \c_0,D^+)$ extends to a $b_q$-sheeted cover $(\tilde U_q,\c _q, D^+ \cup U_q)$.
\end{enumerate}
Set $W\deff U_p\cup U_q$, $W^+ \deff W \cap D^+$, $U_p^+ \deff U_p \cap D^+$ and $U_q^+ \deff U_q \cap D^+$. 
It suffices to prove that there exists an analytic cover over $D^+ \cup W$ which extends $\c_p$ and $\c _q$.
By construction of $U_p$ and $U_q$ there exists $z_*\in U_p\cap U_q\cap D^+$ such that the natural morphisms
\begin{equation*}
\isl_{p,*} : \pi_1(U_p^+\backslash\calr,z_*)\to\pi_1(U_p\backslash\calr,z_*)
\end{equation*} and
\begin{equation*}
\isl_{q,*} : \pi_1(U_q^+\backslash\calr,z_*)\to\pi_1(U_q\backslash\calr,z_*)
\end{equation*}
are surjective, see \refthm{pic_sim}.
Recall the following.
\begin{lem} \label{VK_sum}
Let $X=X_1 \cup X_2$ be an union of path-connected open sets such that $X_1 \cap X_2$ is path-connected.
For every $x_*\in X_1\cap X_2$ the natural homomorphism
\begin{eqnarray}
\pi_1(X_1,x_*)\ast \pi_1(X_2,x_*) \to \pi_1(X,x_*)
\end{eqnarray}
is surjective. Here $\pi_1(X_1,x_*)\ast \pi_1(X_2,x_*)$ denotes the free product of $\pi_1(X_1,x_*)$ and $\pi_1(X_2,x_*)$.
\end{lem}
This is the ``weakest'' part of the proof of Van Kampen Theorem, see \cite{Ha} for more details.
By \lemma{VK_sum} it implies that
every loop $\alpha$ inside $W\backslash\calr$ is homotopic to $\alpha_1^{n_1}\cdots \alpha_p^{n_p}$ where $n_k \in \zz$ and $\alpha_k\subset U_p^+ \backslash \calr$ or  $\alpha_k\subset U_q^+ \backslash \calr$, $k=1,\dots,p$.
It follows that
$\isl_{*} : \pi_1\left(W^+\backslash\calr,z_*\right)\to\pi_1\left(W\backslash\calr,z_*\right)$
and $\isl_{*} : \pi_1\left(D^+\backslash\calr,z_*\right)\to\pi_1\left((D^+ \cup U)\backslash\calr,z_*\right)$ are surjective.
By \refthm{main_thm} there exists a $b_{W}$-sheeted connected analytic cover over $D^= \cup U$ with ramification divisor $\calr$ which extends $ (\tilde D^+, \c_0,D^+)$, where $b_W \leq b_0$.
\lemma{ext_vois} is proved.
\smallskip

\qed

\smallskip

\begin{figure}[h]
\includegraphics[width=0.5\linewidth]{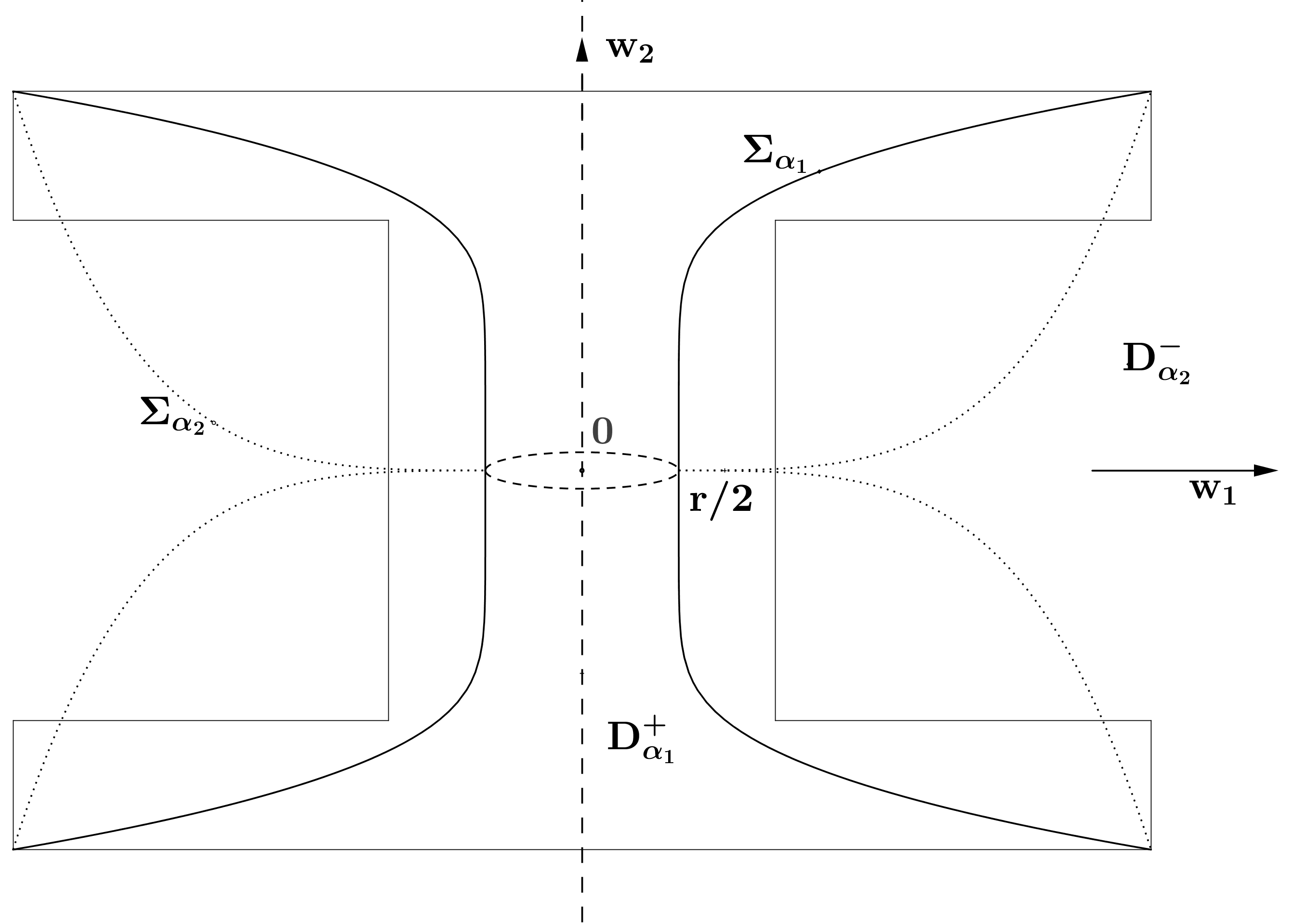}
\caption{For $\alpha_1\gg1$ the hypersurface $\Sigma_{\alpha_1} $
(solid line) belongs to $H^{n,n-q}_{r}$. For $\alpha_2<1$ (dot line)
$\Sigma_{ \alpha_2} $ approaches
$\{w_2=0\}$ when $\alpha_2\searrow 0$. 
Those hypersurfaces are smooth and strictly
$(n-q+1)$-convex outside of the
sphere $\{|w_1|\leq\eps /2, w_2=0\}$ (punctured line).}
\label{gamma-fig}
\end{figure}

\noindent{\slsf Step 3 : }Now we need to exaust $\Delta^n$ by smooth $(n-q+1)$-convex domains starting from $H^{n,n-q}_r$. 
The idea of the following construction is inspired by \S 3 from \cite{Iv2}. 
For $\alpha >0$ consider the  smooth function
\begin{equation} \label{rho-alf}
\rho_{\alpha} (w) = -|w_1|^2 + \frac{r^2}{4} + \left(1-\frac{r^2}{4}\right)|w_2|^{2\alpha}.
\end{equation}
Here $\cc^n\ni w = (w_1;w_2)$ with $w_1=(z_1,\dots,z_{n-q})$, $w_2 = (z_{n-q+1},\dots,z_n)$. 
Set
\begin{equation} \label{d-alfa}
D_{\alpha}^+ = \{ w\in \Delta^n \tq \rho_{\alpha} (w)>0\}, \quad  D_{\alpha}^-\deff \Delta^n\setminus \bar D_{\alpha}^+
\end{equation}
and the hypersurface
\begin{equation} \label{gam-alf}
\Sigma_{\alpha} = \{ w\in \Delta^n:
\rho_{\alpha} (w)=0\}
\end{equation}
separating $D_{\alpha}^+$ from $D_{ \alpha}^-$, see
Figure \ref{gamma-fig}.

\begin{lem}
\label{exhaust-l1}
\begin{enumerate}[label=\roman{*})]
 \item For every $\alpha >0$ the hypersurface $\Sigma_{\alpha}$ is strictly $(n-q+1)$-convex in
$\Delta^n\setminus \{|w_1| \leq \frac{r}{2}, w_2=0\}$. \label{exhaust-l1-i}

\item For $\alpha$ sufficiently big the domain $D^+_{\alpha}$ is contained in $H^{n,n-q}_{r}$ and \label{exhaust-l1-ii}
\begin{equation}
\bigcup_{\alpha >0}D_{\alpha}^+ = \Delta^n\setminus
\left(A^{n-q}_{\frac{r}{2},1}\times\{0\}\right).
\end{equation}
\end{enumerate}

\end{lem}

\proof The Levi form of $\rho_{\eps , \alpha}$ at $w$ is

\begin{equation}
\eqqno(levi-form) 
L_{\rho,{\alpha}}(w)=
\left(%
\begin{array}{cc}
  -2\cdot \adyn_{n-q} & 0  \\
  0 & \left(1-\frac{r^2}{4}\right) H 
  \end{array}%
\right)
\end{equation}
where $H$ is the matrix of the Levi form of $(z_{n-q+1},\dots,z_n)\mapsto \left( |z_{n-q+1}|^2 + \dots +|z_n|^2 \right)^{\alpha}$ which is strictly plurisubharmonic.
This readily implies the assertion \textit{\ref{exhaust-l1-i}} of the Lemma. Assertion \textit{\ref{exhaust-l1-ii}} is left to
the reader.

\smallskip\qed

\smallskip\noindent{\slsf Step 4.} {\it End of the proof.} Applying \lemma{ext_vois} one deduces the set of $\alpha$ such that our cover extends to a neighbourhood of $D_{\alpha}^+$ is a non-empty open closed subset of $]0,+\infty [$. 
By \lemma{exhaust-l1} we extend the cover to $\Delta^n\setminus
\big(A^{n-q}_{\frac{r}{2},1}\times\{0\}\big)$. 
According to Thullen type extension \refthm{ext_thullen} there exists a $b_1$-sheeted connected analytic cover $(\tilde D_1,\c _1,\Delta^n)$ which extends $(\tilde H,\c _0,H_{r}^{n,n-q})$.
By construction that cover verifies conditions \emph{\ref{b'_1<b_1}} and \emph{\ref{b'_1=b_1}} in \refthm{main_thm}. \refthm{ext_hart_b} is proved.

\smallskip


The second assertion is a consequence of \refthm{thm_ext_2}.
Indeed by Theorem of Rothstein (see Corollary 2.19 in \cite{Si}) the ramification divisor $\calr_0$ of the cover uniquely extends to an analytic hypersurface $\calr_1\subset \Delta^n$.
Thus $\calr_1$ is the zero set of a holomorphic function $f\in\calo(\Delta^n)$.
The cover $(\tilde H,\c _0,H_{r}^{n,n-q})$ also extends in the strong sense to a $2$-sheeted connected analytic cover over $\Delta^n$.

\smallskip

\qed

\newsect[EXPLE]{Examples}

\newprg[EXPLE_contre]{An exemple related to Theorem \ref{main_thm}}
We give in this section details about \refexmp{exmp_braid} from Introduction.
\label{cov-incr1} Let $m \geq 4$ be an integer and $X_m$ be the configuration space of $m$ points in $\cc$, \ie $X_m$ is the space of unordered $m$-tuples $(z_1,...,z_m)$ of complex numbers.
It can be identified with $\cc^m$ viewed as the space of monic polynomials of degree $m$.

\begin{equation}
\label{pol}
X^m + \sum_{k=0}^{m-1} w_k X^k
\end{equation}
with coefficients being elementary symmetric functions on $z_1,..,z_m$.
I.e., $w_0 = (-1)^mz_1...z_m, ...,$ $w_{m-1} =- \left(z_1+...+z_m\right)$ or, equivalently $z_1,...,z_m$ are the roots of \eqref{pol}. 
Denote by
\begin{equation}
\label{discr1}
\calr_1 \deff \left\{w\in \cc^m: \discr_X\left(X^m + \sum_{k=0}^{m-1}w_k X^k\right)
=0\right\}
\end{equation}
the discriminant variety, \ie the divisor of such $w$ that the corresponding polynomial has multiply roots. 
Or, in $z$-presentation such $m$-tuples $(z_1,...,z_m)$ where not all $z_j$ are distinct. 
The discriminant variety has equation
\begin{equation}
\eqqno(discr2)
\calr_1 = \left\{z\in X_m: \prod_{i<j}(z_i-z_j)^2=0\right\}.
\end{equation}
Set $D_1\deff X_m$. 
Fix a point $p \deff(\underbrace{1,...,1}_n, 1+\frac{1}{m}, ..., 1 + \frac{m-n}{m})\in \calr_1$, here $3\le n <m$. 
Fix a sufficiently small neighborhood $D_0$ with center at $p$. 
One can suppose $D_0$ is the product of a polydisk $\Delta^n$ with center at $(1,\dots,1)$ in $X_n$ and a polydisk $\Delta^{m-n}$ centered at $(1+\frac{1}{m}, ..., 1 + \frac{m-n}{m})$. 
Set $\calr_0\deff D_0\cap \Delta$. 
Finally fix some $z_*\in D_0\setminus \calr_0$. 
As it is well known $\pi_1(D_1\setminus \calr_1,z_*) = B_m$ the braid group with $m$ strands, see \cite{KT} for details. 
After shrinking $D_0$ one has $\pi_1(D_0\setminus \calr_0,z_*) = B_n$. 
Let $\g _ 0 : B_n \to S_n$ be the natural homomorphism onto the symmetric group of $n$ elements.
It defines a regular $n$-sheeted connected cover $ \c_0 : \tilde D_0'\to D_0\setminus \calr_0$ over $D_0\setminus \calr_0$, see \refrema{corresp_cover}(b).
By \refthm{ext_norm_stein} it can be extended to a $n$-sheeted connected analytic cover over $D_0$ with ramification divisor $\calr_0$.
Suppose that this cover extends to some $n$-sheeted connected analytic cover $\c_1:\tilde D_1\to D_1$ over $D_1$ with divisor $\calr_1$.
Then there should exist a homomorphism $\g_1:B_m\to {S}_n$ making the following diagram commutative.
\begin{equation}\label{diag_n}
\xymatrix{
  B_n \ar@{^{(}->}[r]^{\isl_*} \ar[d]_{\g_0} & B_m\ar[d]^{\g_1} \\
    {S}_n \ar[r]^{\id} & {S}_n
  }
\end{equation}

Here  $\isl_*: B_n \to  B_m$ denotes the homomorphism induced by the natural inclusion $\isl : D_0\backslash \calr_0 \to D_1\backslash \calr_1$.

\begin{prop}
\label{cov-ext1}
For  $(m,n) = (4,3)$ such $\g_1$ doesn't exist.
\end{prop}
\proof Suppose the contrary.
Let $\sigma_1, \sigma_2, \sigma_3$ be the twist generators of $B_4$ and denote by $t_i=\g_1(\sigma_i)\in{S}_3$  their images under $g_1$. 
Then applying $\g_0=\g_1\circ \isl_*$ to $\sigma_1, \sigma_2$ we see that we can suppose that $t_1=(12)$ and $t_2=(23)$. 
Moreover all twists in the braid group are conjugate to each other.
Indeed for every $k=1,\dots,3$ one has $c^k \sigma_1 c^{-k}=\sigma_k$ where $c=\sigma_1\dots\sigma_4$.
It follows that all $t_i$ are conjugate as well. 
Therefore $t_3$ should be a transposition and, moreover, it should commute with $t_1$ and $t_2$. 
But no transposition in ${S}_3$ commutes with $t_1$ and $t_2$. Contradiction.

\smallskip\qed

\begin{rema} \rm\label{cov-incr2}
\begin{enumerate}[label={(\alph{*})}]
\item One can prove more generally that such $\g_1$ doesn't exist for $m>n\ge 3$.

\item Since the diagram 
\begin{equation}\label{nnexmp_diag}
\xymatrix{
  B_n \ar@{^{(}->}[r]^{\isl_*} \ar[d]_{\g_0} & B_m \ar[d]^{\g_1} \\
    {S}_n \ar@{^{(}->}[r]^{\isl} & {S}_m
  }
\end{equation}
is commutative the cover $(\tilde D_0,\c _0,D_0)$ tautologically extends to a $5$-sheeted connected cover over $D_1$.
Here $\isl : S_n \hookrightarrow S_m$ denotes the standard monomorphism of symmetric groups.
\end{enumerate}

\end{rema}

\begin{figure}
    \centering
    \subfigure[Construction of $D_0$ ]{\includegraphics[width=0.3\linewidth]{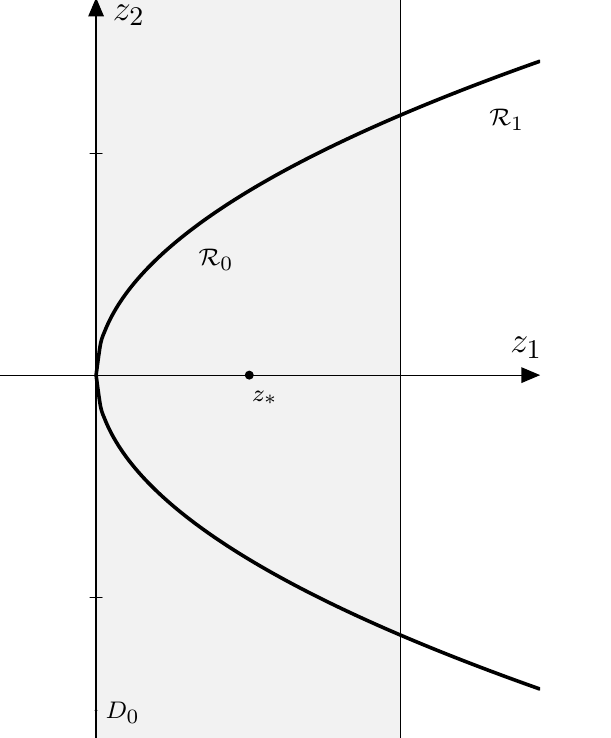}  \label{fig:domains} } \hspace{2cm}
    \subfigure[The cover over $\cc_{z_*}$ ]{ \includegraphics[width=0.5\linewidth]{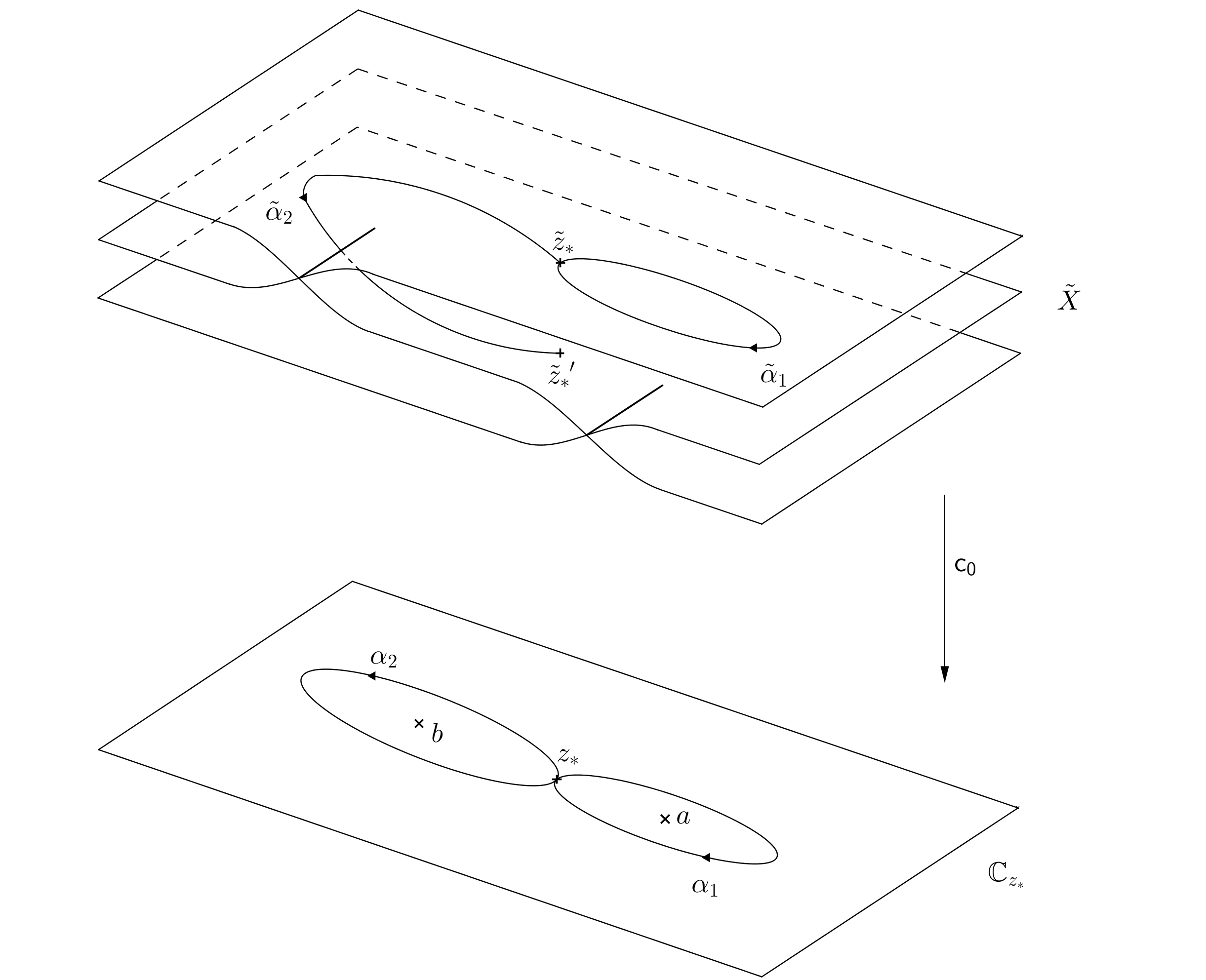} \label{fig:sheets} }
    \caption{}
    \label{fig:cover}
\end{figure}

\newprg[EXPLE_weak]{A non-extendible cover in the strong sense}
We detail here \refexmp{nnex_3} from Introduction.
The idea is inspired by the construction of the cover in Figure \ref{fig:cover}.

Fix $z_*\neq 0$ in $\cc$ and set $\cc_{z_*}\deff \{z_*\}\times\cc$, $D_1 \deff \cc^2$. 
Define $\calr_1\deff \{z_1=z_2^2\}$.
The intersection $\calr_1 \cap \cc_{z_*}$ consists of two points $a$ and $b$.
Let $(\tilde X,\c_0,\cc_{z_*})$ be the $3$-sheeted analytic cover as shown in Figure \ref{fig:sheets}.
There exist two generators $\alpha_1$ and $\alpha_2$ of $\pi_1\left(\cc_{z_*}\backslash\{a,b\},z_*\right)$ such that the lift of $\alpha \deff \alpha_1 \cdot \alpha_2$ at some preimage $\tilde{z}_*$ of $z_*$ is open.
Since $\alpha$ is equivalent to the constant path inside $D_1\backslash\calr_1$ it follows the cover can not be extended to a $3$-sheeted analytic cover over $D_1$.
However one can extend it over $D_0 \deff \Omega \times \cc$ where $\Omega \subset \cc$ is a neighbourhood of $z_*$ such that $\pi_1(D_0 \backslash \calr_1,z_*)= \pi_1(\cc_{z_*}\backslash\{a,b\},z_*)$, see Figure \ref{fig:domains}.

We give here an analytic definition of the cover $(\tilde D_0,\c _0,D_0) $.
Let $f$ be a holomorphic function on  $\Omega \deff \left\{z_1 \in \cc  \tq \left|z_1-1\right|<1  \right\}$ such that $f(z_1)^3=z_1$ and $f(1)=1$.
Set $D_0 \deff \Omega\times \cc$, 
\begin{equation}\label{def_3_cover}
\tilde D_0 \deff \left\{(z,w)\in D_0\times \cc \tq w^3+\frac{f(z_1)}{\sqrt[3]{4}}w+\frac{iz_2}{\sqrt{27}}=0
\right\}
\end{equation}
and $\c_0:\tilde D_0 \to D_0$ induced by the canonical projection $D_0 \times \cc \to D_0$.
Set $\calr_0\deff\calr_1\cap D_0$.
\begin{lem}\label{ex_3_cov}
In the notations as above we have the following.
\begin{enumerate}[label=\emph{(\arabic*)}]
\item  $\left(\tilde D_0,\c _0, D_0\right)$ is a $3$-sheeted connected analytic cover with ramification divisor $\calr_0$; \label{3_cover}

\item  if $(\tilde D_1,\c_1,D_1)$ is a connected extension of $\left(\tilde D_0,\c _0, D_0\right)$ in the weak sense over $D_1$ with ramification divisor ${\calr}_1$ then the number of its sheets is equal to $1$. \label{3_cover_nonext}

\end{enumerate}
\end{lem}

\proof 
{\it Proof of \ref{3_cover}}.
{The space $\ti D_0$ is a complex manifold and the map $\c _0 : \ti D_0 \to D_0$ is holomorphic, proper and surjective.
Recall that the discriminant of the polynomial $w^3 + pw+q$ with respect to $w$ is equal to $-4p^3 - 27 q^2$.
Fix $(z_1,z_2)\in D_0$. $\calr_0$ is the set where the discriminant of  $w^3 + \frac{f(z_1)}{\sqrt[3]{4}}w+\frac{iz_2}{\sqrt{27}}$ with respect to $w$ vanishes.
One deduces that $(\tilde D_0, \c_0, D_0)$ is $3$-sheeted analytic cover with ramification divisor $\calr_0$.}

\noindent{\it Proof of {\ref{3_cover_nonext}}.}
Take $z_* \deff(1,0) \in D_0 \backslash\calr_0$.
The fundamental group of $D_1\setminus\calr_1$ is isomorphic to $\zz$ and $\pi_1(D_0 \setminus\calr_0,z_*)$ is free-generated by the homotopy classes of those loops
\begin{equation} \label{eq_lacets}
\begin{array}{cc}
\alpha_1(t)=\left(1, 1-e^{-i2\pi t} \right)& \\
& \\
\alpha_2(t)=\left(1, -1+e^{i2\pi t} \right)&
\end{array}
\end{equation}
By \refthm{ext_disk_b} there exists a connected extension $(\tilde D_1,\c _1,D_1)$ of $(\tilde D_0,\c _0,D_0)$ over $D_1$ with ramification divisor ${\calr}_1$.
Moreover the number $b_1$ of its sheets can not be larger than $3$.
Let us prove that $b_1=1$ by contradiction.

\smallskip

Suppose that $b_1=3$. $(\tilde D_1,\c _1,D_1)$ is an extension in the strong sense of $(\tilde D_0,\c _0,D_0)$ over $D_1$.
Let $\gamma$ be the loop defined as
\begin{equation}
\gamma(t)=\left(e^{-i2\pi t},0\right).
\end{equation}
One can see that $\alpha_1$ and $\gamma$ are homotopic inside $D_1\backslash\calr_1$. Indeed the map
\[ \fonction{H}{[0,1]^2}{D_1\backslash\calr_1}{(s,t)}{\left[s+(1-s)e^{-i2\pi t},\sqrt{s}\left(1-e^{-i2\pi t}\right)    \right]}
\]
is a continuous deformation between $\alpha_1$ and $\gamma$.
Likewise the loops $\alpha_2$ and $\gamma^{-1}:t\to \gamma(1-t)$ are homotopic inside $D_1\backslash \calr_1$.
It follows that $\alpha\deff \alpha_1 \cdot \alpha_2$ is equivalent to the constant path in $\pi_1(D_1\backslash\calr_1,z_*)$.

\begin{lem}\rm\label{d=3}
There exists a preimage $x_*$ of $z_*$ by $\c _0$ such that the lifted path of $\alpha_1$ at $x_*$ is closed. 
\end{lem}
\proof
Every point of $\calr_0$ has exactly two preimages by $\c _0$.
Indeed for every $(z_1,z_2)\in \calr_0$ there exists $\zeta \in \cc$ such that $z_1=4\zeta ^6$, $f(z_1)=\sqrt[3]{4}\zeta^2$ and $z_2=2\zeta ^3$.
Then 
\begin{equation}\label{exp_order}
w^3 + \frac{f(z_1)}{\sqrt[3]{4}}w+\frac{iz_2}{\sqrt{27}}= w^3 + \zeta ^2w+2i\frac{\zeta^3}{\sqrt{27}}=\left(w-\frac{2i\zeta}{\sqrt{3}}\right)\left(w+\frac{i\zeta}{\sqrt{3}}\right)^2.
\end{equation}
One deduces there exists a preimage of $(1,1)$ which is of order $1$ \ie there exists an open neighbourhood $V$ of $\left(1,1 \right)$ and a connected component $\tilde V$ of $\c^{-1}_0(V)$ such that the restriction $\restriction{\c_0}{\tilde V}:\tilde V \to V$ induced a homeomorphism.
Fix $z \in V\backslash\calr_0$ and let $\lambda$ be a path between $z$ and $z_*$ inside $D_0\backslash\calr_0$.
One can find a loop $\alpha'_1$ starting at $z$ inside $V\backslash\calr_0$ such that $\alpha_1$ is homotopic to $\lambda^{-1}\cdot\alpha'_1\cdot\lambda$.
Let $x\in \tilde V$ be the preimage of $z$ by $\restriction{\c _0}{\tilde V}$ and let $\tilde \alpha'_1$, $\tilde \lambda$ be the lifted paths of $\alpha'_1,\lambda$ by $\c_0$ at $x$.
By construction the lift of $\lambda^{-1}\cdot\alpha'_1\cdot\lambda$ by $\c _0$ at $x_*\deff \tilde{\lambda}(1)$ is closed. \lemma{d=3} is proved.

\qed

\smallskip
One deduces the lift of $\alpha$ by $\c_0$ at $x_*$ is not closed, see Figure \ref{fig:sheets}. 
If not \lemma {d=3} implies that the lifted path of $\alpha_2$ at $x_*$ is closed too.
Since $\pi_1(D_0\backslash\calr_0,z_*)$ is generated by $[\alpha_1]$ and $[\alpha_2]$ it follows that the lift of every loop $\beta$ inside $D_0\backslash\calr_0$ by $\c _0$ at $x_*$ is closed.
Also $x_*$ can not be joined to any other preimage $x'_*$ of $z_*$ inside $\tilde D_0\backslash\tilde{\calr_0}$.
This is not possible because $\tilde{\calr_0}$ does not locally separate $\tilde D_0$.

By construction we can not extend $(\tilde D_0,\c ,D_0)$ in the strong sense to a cover over $D_1$ because the lifted path of $\alpha$ previously defined is not closed in $\tilde D_0$ while $\alpha$ is equivalent to the constant path in $\pi_1(D_1\backslash\calr_1,z_*)$.
It follows that $b_1<3$ and we have a contradiction.

\smallskip
Suppose that $b_1=2$.
We may assume by \refthm{thm_cover_2} that $\tilde{D_1}=\{(z,\xi)\in  D_1\times\cc \tq \xi^2 = z_1-z_2^2  \}$ and $\c_1$ is induced by the natural projection $(z,\xi) \mapsto z$.
There exists a holomorphic map
\[
\fonction{\tilde\isl}{\tilde D_0}{\tilde D_1}{(z,w)}{\left[z,\phi(z,w)\right]}
\] 
where $\phi$ is holomorphic on $\tilde D_0$ such that $\phi(z,w)^2 = z_1-z_2^2$.
Let $\eps >0$ be such that $z_1\deff 4\eps^3 \in \Omega$.
Remark that the projection
\begin{equation*}
\fonction{\pi}{\tilde D_0\cap \{z_1=4 \eps^3\}}{\cc }{(4\eps^3,z_2,w)}{w}
\end{equation*}
induces a biholomorphism.
Thus there exists a holomorphic map
\begin{equation*}
g:w\in \cc\mapsto \phi\left(4\eps ^3,i\sqrt{27}(w^3+\eps w),w\right)
\end{equation*}
such that
$g(w)^2=4\eps^3+27(w^3+\eps w)^2$.
That is not possible because the image of $g$ contains $0$ (take $w^2 = -\eps /3$).
We have a contradiction.
The cover $(\tilde D_1,\c_1,D_1)$ has exactly one sheet \ie $\c _1:\tilde D_1 \to D_1$ is a biholomorphism.
By \refthm{main_thm} if $(\tilde D'_1,\c'_1,D'_1)$ is another extension of $(\tilde D_0,\c_0,D_0)$ the number of its sheets $b'_1$ is equal to one \ie  $(\tilde D'_1,\c'_1,D'_1)$ and $(\tilde D_1,\c_1,D_1)$ are equivalent.
\lemma{ex_3_cov} is proved.

\qed
\begin{rema} \rm
Since $\calr$ does not intersect $\Delta_4\times\left(\Delta_4\backslash\overline{\Delta_{3}}\right)$ the cover $(\tilde D_0,\c _0,D_0)$ can be extended to a $3$-sheeted connected analytic cover over the $1$-convex Hartogs figure 
\[H\deff \left[\Omega\times\Delta_4\right]\cup \Delta_4\times\left(\Delta_4\backslash\overline{\Delta_{3}}\right).\]

\end{rema}
\newprg[EXPLE_galois]{A non-extendible Galois analytic cover}

We give here details about \refexmp{nnex_galois} from Introduction.
Set $\calr_1 \deff \left\{ (z_1,z_2)\in D_1 \tq z_1=z_2^2   \right\}$ and $\calr_0 = \calr_1 \cap D_0$.
The map $\c : \ti D_0 \to D_0$ induces over $D_0 \setminus \calr_0$ a $3$-sheeted regular cover.
Moreover the set $\ti \calr_0 = \c \inv(\calr_0)$ does not locally separate $\ti D_0$.
Indeed let $(z_1,z_2) \in \calr_0$ where $z_1 \in \Omega$.
Without loss of generality we may assume that $z_1=1$.
If $z_2=1$ we denote by $U_r$ the polydisk centered at $(1,1)$ of radius $r>0$.
There exists $r_0 >0$ such that for every $r \in ]0,r_0[$ the preimage set $\c \inv(U_r\setminus \calr_0)$ is homeomorphic to 
\begin{equation*}
\left\{  (W,z) \in \cc \times U\setminus\calr_0 \tq  W^3 =\left( g(z_1)-z_2  \right)^2   \right\}
\end{equation*}
which is connected.
If $z_2=-1$ one deduces in the same way that $\c \inv(U_r\setminus \calr_0)$ is homeomorphic to 
\begin{equation*}
\left\{  (W,z) \in \cc \times U\setminus\calr_0 \tq  W^3 =\left( g(z_1)+z_2  \right)^2   \right\}
\end{equation*}
which is also connected.

By \refthm{struc_normal} the space $\ti D_0$ inherits a unique structure of a normal complex space such that $(\ti D_0,\c _0,D_0)$ becomes a $3$-sheeted analytic cover with ramification divisor $\calr_0$.

By \refthm{main_thm} there exists a connected extension $(\ti D_1, \c _1,D_1)$ of the cover over $D_1 = \cc^2$ with ramification divisor $\calr_1$.
Let us prove by contradiction that the number $b_1$ of its sheets is equal to one.

Let $\alpha_1$ and $\alpha_2$ be the loops defined in \eqref{eq_lacets}.
Then the lifted of $\alpha = \alpha_1 \cdot \alpha_2$ by $\c_0$ based at any preimage of $(1,0)$ is not closed whereas $\alpha$ is homotopic inside $D_1 \setminus \calr_1$ to the constant path.
It follows that the extension is not in the strong sense and $b_1 \leq 2$.

If $b_1=2$ one deduces by \refthm{thm_cover_2} that $\ti D_1$ is equivalent to 
\begin{equation*}
D'_1 = \left\{ (z,\zeta) \in D_1 \times \cc \tq \zeta^2 = z_1-z_2^2  \right\}.
\end{equation*}
and there exists a holomorphic function $\phi : \ti{D}_0 \to \cc$ such that for every $(z,w)\in \ti D_0$ $\phi(z,w)^2=z_1-z_2^2$.
The function 
$h(w)=\phi(w,(\frac{w^3+1}{2})^2, \frac{w^3-1}{2})$ is well-defined holomorphic on the open set $W = \left\{w\in \cc \tq (\frac{w^3+1}{2})^2 \in \Omega \right\}$ which contains $0$..
We have a contradiction because $h(w)^2 = w^3$ for every $w\in W$.

Finally one deduces that $b_1=1$ and every connected extension of $(\ti D_0,\c_0,D_0)$ over $D_1$ is the trivial cover.

\ifx\undefined\bysame
\newcommand{\bysame}{\leavevmode\hbox to3em{\hrulefill}\,}
\fi

\def\entry#1#2#3#4\par{\bibitem[#1]{#1}
{\textsc{#2 }}{\sl{#3} }#4\par\vskip2pt}

\end{document}